\documentclass[reqno,10pt]{amsart}

\usepackage{silence}
\usepackage[english]{babel}
\usepackage{xcolor}
\usepackage{amsthm}
\usepackage{amsmath}
\usepackage{amssymb}
\usepackage{stackengine}
\usepackage{etoolbox}
\usepackage{a4wide}
\usepackage{accents}
\usepackage{palatino}

\usepackage{mathtools}
\mathtoolsset{showonlyrefs,showmanualtags}
\definecolor{citegreen}{rgb}{0,0.6,0}
\definecolor{refred}{rgb}{0.8,0,0}

\usepackage[colorlinks, citecolor=citegreen,linkcolor=refred]{hyperref}
\makeatletter
\AtBeginDocument{%
  \ifdefined\org@@citex
    \patchcmd{\org@@citex}{,\penalty\@m\ }{,\penalty\@m}{}%
      {\PackageWarning{QuasiConvex}{Could not remove spaces between citation numbers}}%
  \else
    \patchcmd{\@citex}{,\penalty\@m\ }{,\penalty\@m}{}%
      {\PackageWarning{QuasiConvex}{Could not remove spaces between citation numbers}}%
  \fi
}
\makeatother
\usepackage{mathrsfs}
\usepackage{esint}

\theoremstyle{plain}
\newtheorem{Teorema}{Theorem}[section]
\newtheorem{corollary}[Teorema]{Corollary}
\newtheorem{proposizione}[Teorema]{Proposition}
\newtheorem{lem}[Teorema]{Lemma}

\newtheorem{ackn}{Acknowledgements\hspace{-.2em}}

\theoremstyle{definition}
\newtheorem{definizione}[Teorema]{Definition}

\theoremstyle{remark}
\newtheorem{oss}[Teorema]{Remark}

\numberwithin{equation}{section}

\def\R{{{\mathbb R}}}

\def\Q{{{\mathbb Q}}}

\def\NN{{{\mathbb N}}}

\def\LLL{{{\mathscr L}}}

\DeclareRobustCommand{\rchi}{{\mathpalette\irchi\relax}}
\newcommand{\irchi}[2]{\raisebox{+1.5pt}{$#1\chi$}}
\newcommand\varempty{\text{\rm{\O}}}

\DeclareMathOperator{\diam}{diam}

\def\XXint#1#2#3{{\setbox0=\hbox{$#1{#2#3}{\int}$} 
\vcenter{\hbox{$#2#3$}}\kern-.5\wd0}}

\newcommand{\weakkk}{\mathrel{\ensurestackMath{\stackon[-.25em]{\xrightharpoonup{\ \ \,}}{\scriptstyle{\,\star\,}}}}} 

\definecolor{MyDarkGreen}{rgb}{0,0.50,0.04}
\definecolor{codegreen}{rgb}{0,0.6,0}
\definecolor{codegray}{rgb}{0.5,0.5,0.5}
\definecolor{codepurple}{rgb}{0.58,0,0.82}
\definecolor{backcolour}{rgb}{0.95,0.95,0.92}

\title{Quasiconvexity and Lower Semicontinuity in the Riemannian Setting}

\date{\today}
 
\author[Aurora Corbisiero]{Aurora Corbisiero}
\address[Aurora Corbisiero]{Department of Mathematics and Applications ``Renato Caccioppoli'', Universit\`a di Napoli Federico II, Italy}
\email{au.corbisiero@studenti.unina.it}

\author[Chiara Leone]{Chiara Leone}
\address[Chiara Leone]{Department of Mathematics and Applications ``Renato Caccioppoli'', Universit\`a di Napoli Federico II, Italy}
\email{chileone@unina.it}

\author[Carlo Mantegazza]{Carlo Mantegazza}
\address[Carlo Mantegazza]{Department of Mathematics and Applications ``Renato Caccioppoli'', Universit\`a di Napoli Federico II, Italy}
\email{carlo.mantegazza@unina.it}

\makeatletter
\patchcmd{\@setaddresses}{\indent}{\noindent}{}{}
\patchcmd{\@setaddresses}{\indent}{\noindent}{}{}
\patchcmd{\@setaddresses}{\indent}{\noindent}{}{}
\makeatother

\subjclass[2020]{Primary 49J45; Secondary 58E30, 26B25}
\keywords{quasiconvexity, lower semicontinuity, Riemannian manifolds, manifold--valued Sobolev maps}
\dedicatory{Dedicated to Nicola Fusco on his 70th birthday}

\begin{document}

\begin{abstract}
We introduce a notion of quasiconvexity for integrands defined on the bundle of linear maps between the tangent spaces of two smooth, connected and complete Riemannian manifolds $(M,g)$ and $(N,h)$, naturally extending the classical Euclidean notion. For continuous integrands, we prove that this Riemannian quasiconvexity is equivalent to the sequential lower semicontinuity of the associated functional
$$
F(u,\Omega)=\int_\Omega f(du)\,d\mu_g
$$
with respect to the weak$^*$ convergence in $W^{1,\infty}(\Omega,N)$, in the intrinsic sense adopted in the paper, for every bounded open subset $\Omega\subseteq M$. We then establish the Riemannian analogues of the Acerbi--Fusco theorems for Carath\'eodory integrands satisfying the corresponding growth conditions, for every $1\leqslant p\leqslant+\infty$, with respect to the weak and weak$^*$ notions of convergence considered here.
\end{abstract}

\maketitle

\section{Introduction}

A fundamental result in the calculus of variations concerning the search for minimizers of variational functionals is the characterization of the sequential lower semicontinuity of integral functionals defined on Sobolev spaces $W^{1,p}(\Omega, \mathbb{R}^m)$, where $\Omega \subseteq \mathbb{R}^n$ is open and bounded. Specifically, for $1 \leqslant p \leqslant +\infty$, lower semicontinuity with respect to the weak convergence (or weak$^*$ convergence if $p=+\infty$) is equivalent to quasiconvexity of the integrand, provided certain growth conditions are satisfied.

Historically, this equivalence was established through several steps: in the pioneering paper~\cite{articoloTonelli}, Tonelli showed that for a twice continuously differentiable function 
$$
f: [a,b]\times\mathbb{R}\times\mathbb{R}\to [0,+\infty)
$$
and for $u\in W^{1,1}(a,b)$, the functional 
\begin{equation*}
F(u)=\int_{a}^bf(x, u(x), u'(x))\,dx
\end{equation*}
is lower semicontinuous with respect to the uniform convergence of curves if and only if the function $f$ is convex in the third variable. This semicontinuity result was later extended to integral functionals on bounded open subsets of $\mathbb{R}^n$ by various authors; Serrin proved in~\cite{Serrin15} that the differentiability assumptions are not actually required. Subsequent improvements of Serrin's theorem were given by De~Giorgi~\cite{degiorgi}, Olech~\cite{olech} and Ioffe~\cite{ioffe}. 

The sufficiency of convexity extends easily to the vectorial case. However, while for $m=1$ the semicontinuity theorem stated above is optimal, in the sense that the convexity of $f(x,s,\xi)$ with respect to $\xi$ is necessary for the lower semicontinuity of $F$, in the vectorial case for 
$n,m>1$, there are functionals (of considerable interest in the theory of nonlinear elasticity, for instance) that are lower semicontinuous without $f$ being convex with respect to the matrix $\xi=(\xi_i^j)$. 

In the case of $m>1$, the condition on $f$ that turned out to be necessary and, with additional assumptions, also sufficient for the lower semicontinuity of integral functionals is quasiconvexity, introduced by Morrey~\cite{Morrey12} in~1952. 

\begin{definizione}\label{quasidef}
A continuous function $f:\mathbb{R}^{{m\times n}}\to \mathbb{R}$ is {\em quasiconvex} if for every ${\xi}\in \mathbb{R}^{{m\times n}} $ and for every bounded open subset $\Omega$ of $\mathbb{R}^{n}$, there holds
\begin{equation}\label{quasiconvexdef}
f({\xi})\leqslant \fint_{\Omega}f({\xi}+D\varphi(x))\,d\LLL^n(x),
\end{equation}
for every function $\varphi\in C^{\infty}_c(\Omega, \mathbb{R}^m)$.\\
A real function $f: \mathbb{R}^n\times \mathbb{R}^m\times \mathbb{R}^{{m\times n}}\to \mathbb{R}$ is {\em quasiconvex in $\xi\in \mathbb{R}^{m\times n}$} if there exists a subset $Z$ of $\mathbb{R}^n$ with $\LLL^n(Z)=0$, such that for every ${x}\in \mathbb{R}^n\setminus Z$ and for every ${s}\in \mathbb{R}^m$ the function $\xi \mapsto f({x},{s},\xi )$ is {\em quasiconvex}.
\end{definizione}

More precisely, Morrey showed that, under some strong regularity assumptions on $f$, its quasiconvexity is equivalent to the sequential weak$^*$ lower semicontinuity in $W^{1,\infty}(\Omega,\mathbb{R}^m)$ of the functional
\begin{equation}\label{Fmorrey}
u\mapsto F(u,\Omega)=\int_{\Omega}f(x,u(x), D u(x))\,d\LLL^n(x).
\end{equation}
Meyers then extended Morrey's result to the setting of $W^{k,p}(\Omega, \mathbb{R}^m)$ spaces in~\cite{Meyers}.

Acerbi and Fusco in~\cite{AcerbiFusco} obtained a significant improvement of this result: they established such an equivalence for Carath\'eodory integrands with appropriate growth conditions in $W^{1,p}(\Omega, \mathbb{R}^m)$, for $1\leqslant p\leqslant +\infty$. We also mention that Marcellini in~\cite{marcellini} presented an alternative proof of this fact. We list below the theorems proved by Acerbi and Fusco.

\begin{Teorema}[Acerbi--Fusco]\label{Fusco1}
Let $f: \mathbb{R}^n\times\mathbb{R}^m\times\mathbb{R}^{{m\times n}}\to \mathbb{R}$ be a Carath\'eodory function satisfying
$$
0\leqslant f(x, s, \xi)\leqslant a(x)+b(s, \xi),
$$
for almost every $x\in \mathbb{R}^n$ and every $s\in \mathbb{R}^m$, $\xi\in \mathbb{R}^{{m\times n}}$, where $a:\R^n\to\R$ is nonnegative and locally integrable and $b:\R^m\times\mathbb{R}^{{m\times n}}\to\R$ is nonnegative and locally bounded.\\
Then, $f$ is quasiconvex in $\xi$ if and only if for every bounded open set $\Omega$ in $\mathbb{R}^n$ the functional $u\mapsto F(u, \Omega)$ is sequentially weak$^*$ lower semicontinuous on $W^{1,\infty}(\Omega, \mathbb{R}^m)$.
\end{Teorema}

\begin{Teorema}[Acerbi--Fusco]\label{Teocap1pfinito}
 Let $1\leqslant p<+\infty$ and $f: \mathbb{R}^n\times\mathbb{R}^m\times\mathbb{R}^{{m\times n}}\to \mathbb{R}$ be a Carath\'eodory function satisfying
 \begin{equation*}
 0\leqslant f(x,s,\xi)\leqslant a(x)+C\big(|s|^p+|\xi|^p\big),
 \end{equation*}
 for almost every $x\in \mathbb{R}^n$ and every $s\in \mathbb{R}^m$, $\xi\in \mathbb{R}^{{m\times n}}$, where $a:\R^n\to\R$ is nonnegative and locally integrable and $C$ is a nonnegative constant.\\
Then, $f$ is quasiconvex in $\xi$ if and only if for every bounded open set $\Omega$ in $\mathbb{R}^n$ the functional $u\mapsto F(u, \Omega)$ is sequentially weakly lower semicontinuous on $W^{1,p}(\Omega, \mathbb{R}^m)$. 
\end{Teorema}

\begin{oss}
Historically, for $p=1$, the sufficiency of quasiconvexity for the sequential weak lower semicontinuity in $W^{1,1}$ goes back to Fusco~\cite[Theorem~3]{Fusco1980}, for Carath\'eodory integrands satisfying a uniform linear growth bound.
\end{oss}

The purpose of this paper is to extend this theory to the Riemannian setting, namely to integral functionals defined on Sobolev maps between two Riemannian manifolds. Such an extension is natural, since vectorial variational problems in which the maps are constrained to take values in a manifold and the domain is itself curved arise in several contexts: in the elasticity of shells and membranes, where the reference configuration is a curved surface; in the theory of liquid crystals and micromagnetics, where the order parameter is a director field or a magnetization taking values in a sphere, possibly on a curved sample; and in nonlinear sigma models of mathematical physics. The model case is the Dirichlet energy, whose critical points are harmonic maps between $(M,g)$ and $(N,h)$, while the functionals considered here provide a natural framework for anisotropic and non--quadratic energies arising in these settings.

We consider two smooth, connected and complete Riemannian manifolds $(M,g)$ and $(N,h)$ and integrands $f$ defined on
\begin{equation*}
\mathscr{L}(TM,TN),
\end{equation*}
the bundle over $M\times N$ whose fiber at $(x,y)$ is the space of linear maps from $T_xM$ to $T_yN$. For continuous integrands, after introducing a Riemannian notion of quasiconvexity extending the classical Euclidean one, we prove that $f$ is quasiconvex in this sense if and only if, for every bounded open subset $\Omega\subseteq M$, the functional
$$
u\mapsto F(u,\Omega)=\int_\Omega f(du)\,d\mu_g
$$
is sequentially lower semicontinuous with respect to the weak$^*$ convergence in $W^{1,\infty}(\Omega,N)$ adopted below, where $\mu_g$ is the Riemannian volume measure of $(M,g)$.

The simultaneous presence of a curved domain and a curved target is the main geometric feature of the problem. The cases in which one of the two manifolds is Euclidean have previously been considered in the literature. A Riemannian domain and a Euclidean target were treated in our previous work~\cite{CLM1}. The complementary case of manifold--valued maps defined on a bounded open subset of $\R^n$ was studied intrinsically by Focardi and Spadaro~\cite{FocardiSpadaro}. An extrinsic
approach, based on relaxation and tangential quasiconvexification, was developed by Dacorogna, Fonseca, Mal\'y and Trivisa~\cite{DFMT}. Maps between two manifolds were also considered by Bernard and Bessi~\cite{BB}, who used Young measures to prove the existence of minimizers for polyconvex integrands in a suitable class of ``Cartesian maps''. A different notion of Riemannian quasiconvexity, in a framework involving differential forms on fiber bundles, was introduced in~\cite{Li}. We discuss the relation between these works and the present one in Section~\ref{quasisec}.

The paper is organized as follows. In Section~\ref{quasisec} we introduce the concept of R--quasiconvexity, prove equivalent formulations and compare it with the usual hierarchy of convexity notions. In Section~\ref{RDefSec} we define the intrinsic weak and weak$^*$ convergence notions used for manifold--valued Sobolev maps and prove the characterization of sequential lower semicontinuity for continuous integrands in $W^{1,\infty}$. 
Finally, in Section~\ref{AFR} we establish the Riemannian analogues of the Acerbi--Fusco theorems for Carath\'eodory integrands satisfying the corresponding growth conditions, for every $1\leqslant p\leqslant+\infty$. For $1<p\leqslant+\infty$, the notion of convergence we will consider coincides with the one of Focardi and Spadaro~\cite{FocardiSpadaro} (extended to our setting), namely strong $L^p$ convergence of the maps together with uniform boundedness of their differentials in $L^p$. For $p=1$, our notion of convergence also requires equiintegrability of the differentials (a property which is implied by the usual weak convergence in $W^{1,1}(\Omega,\R^m)$). In Remark~\ref{FSpone}, we will show by a counterexample that this condition cannot in general be dropped for Carath\'eodory integrands.

We conclude the paper with some directions for further developments of the theory.

\begin{ackn}
The authors wish to thank Giovanni~Alberti, Luigi~Ambrosio, Andrea~Braides, Nicola~Fusco and Marco~Pozzetta for their valuable suggestions, insightful comments and careful reading of preliminary versions of this manuscript.
The last two authors are members of the INdAM--GNAMPA and are partially supported by the PRIN Project 2022E9CF89 ``GEPSO -- Geometric Evolution Problems and Shape Optimization''.
\end{ackn}

\section*{AI usage statement}
During the preparation of this manuscript, the authors used ChatGPT 6.1 Sol (OpenAI), Claude Fable 5.1 (Anthropic) and Gemini Argon (Google) to assist with language editing, LaTeX and bibliographic checks, and the review and organization of mathematical arguments. All AI--generated suggestions were independently checked and revised as appropriate by the authors, who take full responsibility for the content of the manuscript.

\section{Quasiconvexity in the Riemannian context}\label{quasisec}

In this section we introduce the notion of R--quasiconvexity, which extends the classical Euclidean notion of quasiconvexity to the Riemannian setting. In the next section we show that this condition is equivalent to weak$^*$ sequential lower semicontinuity in $W^{1,\infty}$.

\medskip

Let $(M,g)$ and $(N,h)$ be a pair of smooth, connected and complete Riemannian manifolds of dimensions $n$ and $m$, respectively. We define $\mathscr{L}(TM,TN)$ as the space of the linear maps between a tangent space of $M$ and a tangent space of $N$, that is,
$$
\mathscr{L}(TM,TN)=\big\{\alpha:T_xM\to T_yN\,\,\big|\,\,\text{$x\in M$, $y\in N$ and $\alpha$ linear}\big\}.
$$
This space has a natural vector bundle structure over $M\times N$ via the local parametrizations of $M$ and $N$, with fibers given by the vector spaces
$$
\mathscr{L}(T_{x}M,T_{y}N)=\big\{\alpha_x^y:T_xM\to T_yN\,\,\big|\,\,\text{$\alpha_x^y$ linear}\big\}.
$$
The space $\mathscr{L}(TM,TN)$ is clearly locally diffeomorphic to $\R^n\times\R^m\times\R^{m\times n}$. Indeed, choosing a local orthonormal frame $E_1,\dots,E_n$ in a neighborhood $U\subseteq M$ (diffeomorphic to $\R^n$) of a point $x_0\in M$ and another one $F_1,\dots,F_m$ in a neighborhood $V\subseteq N$ (diffeomorphic to $\R^m$) of $y_0\in N$, in order to locally ``trivialize'' the vector bundle, we define the map $I:\mathscr{L}(TU,TV)\to\R^{m\times n}$ as
\begin{equation}\label{mapI}
\alpha\mapsto I\alpha=\big(\alpha^j(E_i)\big)=\big(h(\alpha(E_i),F_j)\big)\in \R^{m\times n},
\end{equation}
for every $\alpha\in\mathscr{L}(TU,TV)$, hence the map sending any $\alpha_x^y:T_xM\to T_yN$ to 
$$
(x,y,I\alpha_x^y)=\big(x,y,h(\alpha_x^y(E_i),F_j)\big)\in U\times V\times\R^{m\times n}\approx\R^n\times\R^m\times\R^{m\times n}
$$
is a diffeomorphism.

Then, we have a distance $\delta$ on $\mathscr{L}(TU,TV)$ defined as follows: 
\begin{equation}\label{distL}
\delta(\alpha,\beta)=d^M(x,x')+d^N(y,y')+\Vert I\alpha-I\beta\Vert=d^M(x,x')+d^N(y,y')+\sqrt{\sum_{i=1}^n\sum_{j=1}^m \big|\alpha^j(E_i)-\beta^j(E_i)\big|^2\,}
\end{equation}
for every $\alpha,\beta\in\mathscr{L}(TU,TV)$ such that $\alpha\in\mathscr{L}(T_{x}M,T_{y}N)$ and $\beta\in\mathscr{L}(T_{x'}M,T_{y'}N)$, where $d^M$ and $d^N$ are respectively the Riemannian distances on $M$ and $N$.\\
It is then easy to see that such a distance, restricted to every fiber $\mathscr{L}(T_{x}M,T_{y}N)$, coincides with the one associated with the ``quadratic'' norm (equivalent to the operator norm) induced by the metric tensors $g_x$ and $h_y$ of $M$ and $N$.

Clearly, if we have a differentiable function $u: M\to N$, its differential $du[x]:T_xM\to T_{u(x)}N$ belongs to $\mathscr{L}(T_xM,T_{u(x)}N)\subseteq\mathscr{L}(TM,TN)$, for every $x\in M$.
\medskip

{\em We choose to adopt the notation $du[x]$ for the differential of a map $u: M\to N$ at a point $x\in M$, for the sake of clarity in the computations that follow. Moreover, again for notational clarity, we decided to deal throughout the paper with the differentials of the functions instead of their ``gradients''. The two points of view are clearly equivalent.}

\medskip

In the notation above, for every $x\in M$ and $y\in N$, the vector spaces $T_xM$ and $T_yN$ are respectively isometric to $\R^n$ and $\R^m$, via the orthonormal bases $e_1,\dots,e_n$ and $f_1,\dots,f_m$ given by the local orthonormal frames $E_1,\dots,E_n$ and $F_1,\dots,F_m$ at $x$ and $y$ (in a local ``trivialization'' of the bundle $\mathscr{L}(TM,TN)$ around $(x,y)$).\\
With this identification, for $r>0$, we define the open cube
\begin{equation*}
\widetilde{Q}^r_x=r\bigg(-\frac{1}{2}, \frac{1}{2} \bigg)^n\subseteq T_{x}M\approx\R^n
\end{equation*}
and we will always assume that $r$ is small enough that $\widetilde{Q}^r_x$ is contained in a ball of $T_xM$ centered at the origin and with radius smaller than the injectivity radius of $M$ at $x$. This clearly implies that the exponential map
$$
\operatorname{exp}_{x}: T_{x}M\to M
$$
gives a diffeomorphism between $\widetilde{Q}^r_x$ and its image
\begin{equation*}
 Q^r_{x}=\operatorname{exp}_{x}\big(\widetilde{Q}^r_x\big)\subseteq M.
\end{equation*}
In an entirely analogous way, for every $y\in N$, we define the open cube 
$$
\widetilde{P}^s_y=s\bigg(-\frac{1}{2}, \frac{1}{2} \bigg)^m\subseteq T_yN\approx\R^m
$$
and its image $P^s_{y}=\operatorname{exp}_{y}\big(\widetilde{P}^s_y\big)\subseteq N$ via the exponential map $\operatorname{exp}_{y}:T_yN\to N$. We underline, in particular, that the exponential maps $\exp_x$ and $\exp_y$ have inverses defined on ${Q}^r_{x}$ and ${P}^s_{y}$, respectively.

\begin{proposizione}\label{RiemannianQC-equiv}
Let $(M,g)$ and $(N,h)$ be a pair of smooth, connected and complete Riemannian manifolds of dimensions $n$ and $m$, respectively, and let $\mu_g$ be the Riemannian volume measure of $(M,g)$. 
Let $f: \mathscr{L}(TM,TN)\to \mathbb{R}$ be a real function and let $x_0\in M$, $y_0\in N$ be such that the ``restricted'' map
$$
f_{x_0}^{y_0}=f|_{\mathscr{L}(T_{x_0}M,T_{y_0}N)}:\mathscr{L}(T_{x_0}M,T_{y_0}N)\to\R
$$
is continuous. Then, for every $s>0$ small enough that $\exp_{y_0}$ is a diffeomorphism on $\widetilde P_{y_0}^s$, the following three conditions are equivalent:
\begin{itemize}
\item for every $\alpha_{x_0}^{y_0}\in\mathscr{L}(T_{x_0}M,T_{y_0}N)$, every $r>0$ small enough and every 
$\varphi\in C^\infty(Q_{x_0}^r,P_{y_0}^s)$ which is equal to $y_0$ outside a compact subset 
of $Q_{x_0}^r$, there holds
\begin{equation}\label{quasidefR}
f(\alpha_{x_0}^{y_0})\leqslant\fint_{Q_{x_0}^r} f\big(\alpha_{x_0}^{y_0}
+d\exp_{y_0}^{-1}[\varphi(x)]\circ d \varphi [x]\circ d\exp_{x_0}[\exp_{x_0}^{-1}(x)]\big)\,d\mu_g(x)+o(1),
\end{equation}
where $o(1)$ is a function which goes to zero as $r\to 0$ and depends, in a monotonically nondecreasing way, only on the $L^{\infty}$ norm of $d\varphi$ with respect to local orthonormal frames on $Q_{x_0}^r$ and $P_{y_0}^s$, for fixed $x_0$, $y_0$, $s$ small enough and $\alpha_{x_0}^{y_0}$.

\smallskip

\item for every $\alpha_{x_0}^{y_0}\in\mathscr{L}(T_{x_0}M,T_{y_0}N)$, every $r>0$ small 
enough and every $\varphi\in C^\infty(Q_{x_0}^r,P_{y_0}^s)$ which is equal to $y_0$ outside a compact subset of $Q_{x_0}^r$, there holds
\begin{equation}\label{quasidefR2}
f(\alpha_{x_0}^{y_0})\leqslant\frac{1}{r^n}\int_{Q_{x_0}^r} f\big(\alpha_{x_0}^{y_0}
+d\exp_{y_0}^{-1}[\varphi(x)]\circ d \varphi [x]\circ d\exp_{x_0}[\exp_{x_0}^{-1}(x)]\big)\,J_{\exp_{x_0}^{-1}}(x)\,d\mu_g(x),
\end{equation}
where $J_{\exp_{x_0}^{-1}}(x)$ denotes the Jacobian of the map $\exp_{x_0}^{-1}$ at the point $x\in Q_{x_0}^r$.

\smallskip

\item the ``restricted'' map $f_{x_0}^{y_0}$ is quasiconvex in the usual sense: for every 
$\alpha_{x_0}^{y_0}\in\mathscr{L}(T_{x_0}M,T_{y_0}N)$, every $r>0$ and every 
$\psi\in C^\infty_c\big(\widetilde{Q}^r_{x_0},T_{y_0}N\big)$, there holds
\begin{equation}\label{quasidefR3}
f(\alpha_{x_0}^{y_0})\leqslant\fint_{\widetilde{Q}^r_{x_0}} f\big(\alpha_{x_0}^{y_0}+d\psi[z]\big)\,d\LLL^n(z).
\end{equation}
\end{itemize}
\end{proposizione}
\begin{proof}
The second condition implies the first one, since $\mu_g(Q_{x_0}^r)/r^n\to1$ as $r\to0$ and the Jacobian $J_{\exp_{x_0}^{-1}}(x)$ of the inverse exponential map converges uniformly to $1$ (the differential of the exponential map is the identity at the origin of $T_{x_0}M$), as $x\in Q_{x_0}^r$ and $r\to 0$. Moreover, if the $L^{\infty}$ norm of $d\varphi$ is bounded, the integrand $f\big(\alpha_{x_0}^{y_0}+d\exp_{y_0}^{-1}[\varphi(x)]\circ d \varphi [x]\circ d\exp_{x_0}[\exp_{x_0}^{-1}(x)]\big)$ is also bounded.\\
The third condition implies the second one, by setting $\psi=\exp_{y_0}^{-1}\circ\varphi\circ\exp_{x_0}$ and changing variables $z=\exp^{-1}_{x_0}(x)$. Hence, we only need to show that the first condition implies the third one.\\
Assume the first condition holds; then, by what we said above about $J_{\exp_{x_0}^{-1}}$, we also have
\begin{equation}
f(\alpha_{x_0}^{y_0})\leqslant\frac{1}{r^n}\int_{Q_{x_0}^r} f\big(\alpha_{x_0}^{y_0}
+d\exp_{y_0}^{-1}[\varphi(x)]\circ d \varphi [x]\circ d\exp_{x_0}[\exp_{x_0}^{-1}(x)]\big)\,J_{\exp_{x_0}^{-1}}(x)\,d\mu_g(x)+o(1),
\end{equation}
hence, 
\begin{equation}\label{eqcar8881}
f(\alpha_{x_0}^{y_0})\leqslant\fint_{\widetilde{Q}^r_{x_0}} f\big(\alpha_{x_0}^{y_0}+d\psi[z]\big)\,d\LLL^n(z)+o(1),
\end{equation}
where $\psi=\exp_{y_0}^{-1}\circ\,\varphi\circ\,\exp_{x_0}$. Thus, this inequality holds for every $\psi\in C^\infty_c\big(\widetilde{Q}^r_{x_0},\widetilde{P}^s_{y_0}\big)$, $r>0$ small enough, where $o(1)$ is as in the statement, that is, it goes to zero as $r\to 0$ and depends only on the $L^{\infty}$ norm of $d\psi$.\\
Let us fix $r_0>0$ and $\psi\in C^\infty_c\big(\widetilde{Q}^{r_0}_{x_0},T_{y_0}N\big)$. For any $\rho > 0$, we define the rescaled function $\psi_\rho$ on $\widetilde{Q}_{x_0}^{\rho r_0}$ as
\begin{equation}
\psi_\rho(z)= \rho\,\psi\left({z}/{\rho}\right),
\end{equation}
which clearly satisfies
$$
d \psi_\rho[z] = d \psi[{z}/{\rho}],\qquad\text{ hence, }\qquad\|d \psi_\rho\|_{L^\infty} = \|d \psi\|_{L^\infty}.
$$
For $\rho>0$ small enough, the image of $\psi_\rho$ is contained in $\widetilde P_{y_0}^s$; then, by formula~\eqref{eqcar8881} applied on the cube $\widetilde{Q}^{\rho r_0}_{x_0}$ and performing the change of variables $w = z/\rho$ (so that $d\LLL^n(z)= \rho^n d\LLL^n(w)$), we have
\begin{align*}
f(\alpha_{x_0}^{y_0})\leqslant&\,\fint_{\widetilde{Q}^{\rho r_0}_{x_0}} f\big(\alpha_{x_0}^{y_0}+d\psi_\rho[z]\big)\,d\LLL^n(z)+o(1)=\fint_{\widetilde{Q}_{x_0}^{\rho r_0}} f\big(\alpha_{x_0}^{y_0}+d\psi[z/\rho]\big)\,d\LLL^n(z)+o(1)\\
=&\,\frac{1}{r_0^n\rho^n} \int_{\widetilde{Q}_{x_0}^{r_0}} f\big(\alpha_{x_0}^{y_0}+d\psi[w]\big)\,\rho^nd\LLL^n(w)+o(1)
=\fint_{\widetilde{Q}_{x_0}^{r_0}} f\big(\alpha_{x_0}^{y_0}+d\psi[w]\big)\,d\LLL^n(w)+o(1).
\end{align*}
Hence, sending $\rho\to 0$, as the $L^\infty$ norms of all $d\psi_\rho$ are uniformly bounded, the ``error term'' $o(1)$ disappears and we conclude
$$
f(\alpha_{x_0}^{y_0})\leqslant\fint_{\widetilde{Q}_{x_0}^{r_0}} f\big(\alpha_{x_0}^{y_0}+d\psi[w]\big)\,d\LLL^n(w),
$$
that is, the third condition is satisfied and we are done.
\end{proof}

Then, we are ready to give a definition of quasiconvexity in the Riemannian setting.

\begin{definizione}\label{RiemannianQC2}
Let $(M,g)$ and $(N,h)$ be a pair of smooth, connected and complete Riemannian manifolds of dimensions $n$ and $m$, respectively, and let $\mu_g$ be the Riemannian volume measure of $(M,g)$. We call a real function $f: \mathscr{L}(TM,TN)\to \mathbb{R}$ {\em Riemannian--quasiconvex} or simply {\em R--quasiconvex} if for $\mu_g$--almost~every $x_0\in M$ and every $y_0\in N$, the ``restricted'' map $f_{x_0}^{y_0}$ is continuous and $f$ satisfies one of the equivalent conditions in the previous proposition.
\end{definizione}

The main obstacle in generalizing to the Riemannian context the usual definition of quasiconvexity in the Euclidean setting is due to the fact that we cannot identify all the tangent spaces of a manifold as we can clearly do in $\R^n$, hence two differentials at two different points cannot be added together. Thus, in the definition above (referring to the first condition in Proposition~\ref{RiemannianQC-equiv}) we needed to morally ``carry'' the differential of the ``perturbation'' $\varphi$ at every point $x\in Q_{x_0}^r$ back to the fiber $\mathscr{L}(T_{x_0}M,T_{y_0}N)$, via the differentials of the exponential maps at $x_0\in M$ and $y_0 \in N$ (exponential maps whose differentials in the Euclidean case are simply the identity, under the usual identification of $\R^n$ with all its tangent spaces -- see also the discussion below). This forces the introduction of an ``error term'' in the ``quasiconvexity inequality''~\eqref{quasidefR} that one reasonably expects to go to zero as we get closer and closer to the point $x_0$, since the Riemannian volume is only asymptotically Euclidean or, equivalently, because the Jacobian $J_{\exp_{x_0}^{-1}}$ in the
second formulation tends to one.

\begin{oss}
We underline that the presence of the ``error term'' $o(1)$ in formula~\eqref{quasidefR} is only due to the use of the exponential map in dealing with the non--flatness of the domain $(M,g)$, not of the target space $(N,h)$, as shown below by the Focardi--Spadaro case, where the domain is flat. Indeed, by looking at the second equivalent condition of R--quasiconvexity in Proposition~\ref{RiemannianQC-equiv}, if the Jacobian $J_{\exp_{x_0}^{-1}}(x)$ of the inverse exponential map is identically one, as in the Euclidean case, then the ``error term'' $o(1)$ is not necessary.
\end{oss}

\begin{oss}\label{indiposs}
It might appear that the conditions in Proposition~\ref{RiemannianQC-equiv} (hence Definition~\ref{RiemannianQC2}) depend on the choice of the orthonormal bases that we use to identify the tangent spaces with $\R^n$ and $\R^m$, in order to define the cubes $Q_{x_0}^r$ and $P_{y_0}^s$, since changing such bases clearly modifies such cubes. However, by the third equivalent condition in Proposition~\ref{RiemannianQC-equiv} and the usual domain--independence of quasiconvexity, this dependence is only apparent. The same argument shows that, instead of cubes, we can work with balls in the tangent spaces and their images via the exponential maps (hence, geodesic/metric balls of the manifolds). We opted for cubes since this choice simplifies the proofs of the theorems of the next section.
\end{oss}

We briefly comment on the three equivalent definitions of R--quasiconvexity given by Proposition~\ref{RiemannianQC-equiv}. The first definition, which is the most similar to the ``classical'' one, being given in terms of an integration on the domain of the functions, can be seen as a sort of ``asymptotic quasiconvexity''. Moreover, the presence of the ``error term'' $o(1)$ in inequality~\eqref{quasidefR} is useful in the approximation arguments employed in the proof of the characterization of lower semicontinuous functionals in the next section.\\
The third one, which deals with an integration on the tangent space of $M$, is useful for comparing R--quasiconvexity with other notions of convexity, as discussed below. Since Proposition~\ref{RiemannianQC-equiv} is stated at a fixed pair $(x_0,y_0)$, Definition~\ref{RiemannianQC2} can then be rephrased as follows: for $\mu_g$--almost~every $x_0\in M$ and every $y_0\in N$, the map $f_{x_0}^{y_0}$ is quasiconvex in the usual sense.\\
Finally, the second equivalent condition in Proposition~\ref{RiemannianQC-equiv} provides a bridge connecting the first formulation (focused on the manifold) and the third one (privileging the
tangent space).

\medskip

If $(M,g)$ is $\R^n$ with its standard metric, for every point $x\in\R^n$ we have a standard identification $T^*_x\R^n\approx T_x{\R^n}\approx\R^n$ and similarly for $(N,h)=\R^m$. Hence, a function $f:\mathscr{L}(T\R^n,T\mathbb{R}^m)\to \mathbb{R}$ can be put in one--to--one correspondence with a function $\widetilde{f}:\R^n\times\R^m\times\R^{m\times n}\to\R$ as follows: if $A=(A_i^j)\in \R^{m\times n}$, for $i\in\{1,\dots,n\}$ and $j\in\{1,\dots,m\}$, then, 
\begin{equation}\label{eq560}
\widetilde{f}({x_0},y_0,A)=f(\alpha_{x_0}^{y_0}),
\end{equation}
where $\alpha_{x_0}^{y_0}:T_{x_0}\R^n\to T_{y_0}\R^m$ is given by setting $\alpha^{y_0}_{x_0}(v)$ equal to the vector of $T_{y_0}\R^m\approx\R^m$ with components
$$
\alpha^{y_0}_{x_0}(v)^j=\sum_{i=1}^n A_i^j v^i,
$$
for every vector $v=(v^1,\dots,v^n)\in T_{x_0}\R^n\approx\R^n$ and viceversa: if $\alpha_{x_0}^{y_0}$ is operating as in this formula, then $f$ is defined by equality~\eqref{eq560}.\\
Then, with this one--to--one correspondence, our ``Riemannian'' definition of quasiconvexity in the case of the Euclidean space (the first condition in Proposition~\ref{RiemannianQC-equiv}) is equivalent to 
\begin{equation}\label{eq562}
\widetilde{f}\big(x_0,y_0,A\big)\leqslant\fint_{Q_{x_0}^r} \widetilde{f}\big(x_0,y_0,A+d\varphi[x]\big)\,d\LLL^n(x)+o(1),
\end{equation}
for $\LLL^n$--almost~every $x_0\in\R^n$, every $y_0\in\R^m$, $A\in\R^{m\times n}$ and $\varphi\in C^\infty_c(Q_{x_0}^r,\mathbb{R}^m)$, since, under the identification $T^*_x\R^n\approx T_x{\R^n}\approx\R^n$, the differentials of all the exponential maps are the identity. Moreover, inequality~\eqref{eq562} is clearly satisfied by a quasiconvex function $\widetilde{f}:\R^n\times\R^m\times\R^{m\times n}\to\R$, actually without the ``error term'' $o(1)$, according to the usual Definition~\ref{quasidef} of quasiconvexity.\\
Hence, a ``classical'' quasiconvex function $\widetilde{f}$ in the Euclidean setting is also R--quasiconvex in the sense of Definition~\ref{RiemannianQC2} (with the above identification of $\widetilde{f}$ with $f$). Conversely, we can directly resort to the equivalent definition of R--quasiconvexity given by the third condition in Proposition~\ref{RiemannianQC-equiv}. Indeed, with the identifications
$T^*_{x_0}\mathbb R^n\approx T_{x_0}\mathbb R^n\approx\mathbb R^n$ and $T^*_{y_0}\mathbb R^m\approx T_{y_0}\mathbb R^m\approx\mathbb R^m$, this condition becomes precisely Definition~\ref{quasidef}.\\
Thus, our definition is an extension to Riemannian manifolds of the usual definition of quasiconvexity in the standard Euclidean context.

\medskip

In our previous work~\cite{CLM1}, we addressed the simpler case when $(N,h)$ is $\R^m$ with its standard metric, while only the domain of the functions is a Riemannian manifold, obtaining the characterization proved in the next section, in the case of continuous integrands and $W^{1,\infty}$ maps. Then, in the definition of R--quasiconvexity, one does not need the composition with $d\exp_{y_0}^{-1}[\varphi(x)]$ appearing in formula~\eqref{quasidefR}, since it is possible to identify all the tangent spaces of $\R^m$ with $\R^m$ itself. On the other hand, Focardi and Spadaro in~\cite{FocardiSpadaro} gave a definition of quasiconvexity in the ``complementary'' case of Sobolev maps $u:\Omega\to N$, with $\Omega$ a bounded open subset of $\mathbb{R}^n$ and $(N,h)$ a Riemannian manifold, which is actually the third condition in Proposition~\ref{RiemannianQC-equiv} (hence, no ``error term'' is present in their definition) and, in the $W^{1,\infty}$ case, proved the analogous equivalence with the weak$^*$ lower semicontinuity of the functionals, with respect to the intrinsic notion of convergence adopted there.

The same constrained problem had been considered before, from an ``extrinsic'' point of view, by Dacorogna, Fonseca, Mal\'y and Trivisa in~\cite{DFMT}: embedding the target manifold $N$ into a Euclidean space $\R^d$, they obtained an integral representation for the relaxation of the energy, whose density is the ``tangential quasiconvexification'' of the integrand, that is, the quasiconvexification obtained by using only perturbations with values in the tangent space of $N$. The case of maps between two manifolds was instead considered by Bernard and Bessi in~\cite{BB}, who proved, by means of Young measures, the existence of a minimizer for polyconvex integrands in a suitable class of ``Cartesian maps''.

Finally, we mention that another notion of Riemannian quasiconvexity, in a
different framework involving differential forms on fiber bundles, was introduced
in~\cite{Li}.

\begin{oss}
We observe, as Focardi and Spadaro do in~\cite{FocardiSpadaro}, that the third equivalent definition of R--quasiconvexity in Proposition~\ref{RiemannianQC-equiv} is related only to the tangent spaces of the two manifolds $M$ and $N$. Hence, it could possibly be extended to the setting of maps between metric spaces, provided the latter possess a suitable notion of tangent spaces (and possibly of an analogue of the exponential map).
\end{oss}

In the Euclidean setting, it is well known that for a function $f:\R^{m\times n}\to\R$, the following ``hierarchy'' holds (see~\cite[Theorem~5.3]{dacorogna}):
$$
\text{Convex} \quad\implies\quad \text{Polyconvex} \quad\implies\quad \text{Quasiconvex} \quad\implies\quad \text{Rank--one convex}
$$
Then, we can ask ourselves if these implications are still valid in our Riemannian context, after giving ``natural'' generalizations of these notions of ``convexity''. The following definitions are an adaptation to the Riemannian setting of the corresponding Euclidean definitions, which can be found in~\cite[Chapter~5]{dacorogna}. It is easy to check that they all coincide with the usual ones when $(M,g)$ and $(N,h)$ are respectively $\R^n$ and $\R^m$ with their canonical metrics.

\begin{definizione}\label{convdef}
Let $(M,g)$ and $(N,h)$ be a pair of smooth, connected and complete Riemannian manifolds of dimensions $n$ and $m$, respectively, and let $\mu_g$ be the Riemannian volume measure of $(M,g)$. 
We call a function $f: \mathscr{L}(TM,TN)\to \mathbb{R}$ {\em R--convex}, {\em R--polyconvex} or {\em R--rank--one convex} if for $\mu_g$--almost~every $x_0\in M$ and every $y_0\in N$, the ``restricted'' map
$$
f_{x_0}^{y_0}=f|_{\mathscr{L}(T_{x_0}M,T_{y_0}N)}:\mathscr{L}(T_{x_0}M,T_{y_0}N)\to\R
$$
is respectively convex, polyconvex or rank--one convex in the usual sense.
\end{definizione}

Then, recalling from the discussion above the fact that the third condition in Proposition~\ref{RiemannianQC-equiv} is equivalent to the usual quasiconvexity for the ``restricted'' map $f_{x_0}^{y_0}$, we can immediately conclude that for a function $f: \mathscr{L}(TM,TN)\to \mathbb{R}$ there holds:
$$
\text{R--convex} \quad\implies\quad \text{R--polyconvex} \quad\implies\quad \text{R--quasiconvex} \quad\implies\quad \text{R--rank--one convex}
$$

\begin{corollary}
If $f:\mathscr{L} (TM,TN)\to \mathbb{R}$ is R--quasiconvex, then, for $\mu_g$--almost~every $x_0\in M$ and every $y_0\in N$, the ``restricted'' map $f_{x_0}^{y_0}$ is locally Lipschitz.
\end{corollary}
\begin{proof}
The claim holds for any R--rank--one convex function, hence, by the previous discussion, also for an R--quasiconvex function. Indeed, rank--one convexity implies separate convexity; therefore, by~\cite[Theorem~2.31]{dacorogna}, the map $f_{x_0}^{y_0}$ is locally Lipschitz.
\end{proof}

\begin{oss}
We also observe that if $f$ is R--quasiconvex, then, for $\mu_g$--almost~every $x_0\in M$ and every $y_0\in N$ such that $f_{x_0}^{y_0}$ is of class $C^2$, the map $f_{x_0}^{y_0}$ satisfies the Legendre--Hadamard condition (relevant for the regularity of the minimizers, see for instance~\cite{giusti}
and~\cite{GiaquintaMartinazzi}): for every $\alpha_{x_0}^{y_0}\in {\mathscr{L}(T_{x_0}M,T_{y_0}N)}$, $\eta\in T_{x_0}^*M$ and $\lambda\in T_{y_0}N$, expressing everything in terms of coordinates given by some orthonormal bases of $T_{x_0}M$ and $T_{y_0}N$, there holds
\begin{equation}\label{eq:LegendreHadamard}
\sum_{i,j=1}^n \sum_{k,\ell=1}^m\frac{\partial^2 f_{x_0}^{y_0}}{\partial\xi_{i}^k \partial\xi_{j}^\ell}\big(\alpha_{x_0}^{y_0}\big) \,\eta_i \eta_j\lambda^k \lambda^\ell \geqslant 0.
\end{equation}
\end{oss}

\section{R--quasiconvexity and weak$^*$ semicontinuity in $W^{1,\infty}$}\label{RDefSec}

In this section we first consider the case of $W^{1,\infty}$ maps and continuous integrands, in order to focus on the geometric aspects related to the Riemannian structure of both the domain and the target, without introducing at the same time the additional technical difficulties arising from the $W^{1,p}$ setting and from a merely measurable dependence of the integrand on the point $x\in M$. This ``continuous case'' will also be used directly in Section~\ref{AFR}: after a Scorza~Dragoni localization and freezing the dependence on the base point, the sufficiency proof for Carath\'eodory integrands in $W^{1,\infty}$ reduces locally to Theorem~\ref{QCR->sci}. The more general case of $W^{1,p}$ maps and Carath\'eodory integrands is treated in the next section. We prove, in our Riemannian setting, the analogues of some special cases of Theorem~\ref{Fusco1} of Acerbi and Fusco, following the line of~\cite{Ambrosio} (we mention that the main argument in the next Theorem~\ref{QCR->sci} is due to Fonseca--M\"uller~\cite{FonsecaMuller}). 

Notice that, since the base of the bundle $\mathscr{L}(TM,TN)$ is $M\times N$, the integrand $f(du[x])$, with $du[x]\in\mathscr{L}(T_xM,T_{u(x)}N)$, depends on $x$, on $u(x)$ and on $du[x]$, hence in the Euclidean case it corresponds to a continuous integrand $f(x,u,Du)$.

\medskip

We recall the Sobolev spaces of maps with values in a Riemannian manifold and we introduce the notions of weak and weak$^*$ convergence used throughout the paper.

Following~\cite{FocardiSpadaro}, for every bounded open set $\Omega\subseteq M$ and $1\leqslant p\leqslant+\infty$, we define $W^{1,p}(\Omega,N)$ intrinsically as follows. A measurable map $u:\Omega\to N$ belongs to $W^{1,p}(\Omega,N)$ if there exists a function $H\in L^p(\Omega)$ such that, for every $y_0\in N$, the function
$$
x\mapsto d^N\big(u(x),y_0\big)
$$
belongs to $W^{1,p}(\Omega)$ and
$$
\big|d\big(d^N(u(\cdot),y_0)\big)[x]\big|_{g_x}\leqslant H(x)
$$
for $\mu_g$--almost~every $x\in\Omega$. The role of the function $H$ is to require a bound on the differentials of the functions $d^N(u(\cdot),y_0)$ which is uniform with respect to $y_0\in N$: each of them only ``sees'' the component of the differential of $u$ in the direction of the geodesic from $y_0$ to $u(x)$, hence asking that all of them belong to $W^{1,p}(\Omega)$, with no common bound, would give no control on the map $u$ itself. As shown in~\cite{FocardiSpadaro}, the optimal choice of $H$ is the function
$$
x\mapsto\sup_{i\in\NN}\big|d\big(d^N(u(\cdot),y_i)\big)[x]\big|_{g_x},
$$
for any countable dense subset $\{y_i\}_{i\in\NN}$ of $N$. As in~\cite{FocardiSpadaro}, every $u\in W^{1,p}(\Omega,N)$ is approximately differentiable for $\mu_g$--almost~every $x\in\Omega$ and its approximate differential defines a linear map
$$
du[x]\in\mathscr L(T_xM,T_{u(x)}N).
$$
We denote by $|du[x]|$ the norm of $du[x]$ induced by $g_x$ and $h_{u(x)}$. The proof of the almost everywhere existence of the approximate differential is local and the argument of~\cite{FocardiSpadaro} applies in local coordinates on $M$. Moreover, again as in~\cite{FocardiSpadaro}, the norm $|du|$ is bounded by $H$, up to a constant depending only on $n$ and $m$; in particular, $|du|\in L^p(\Omega)$. As recalled in~\cite{FocardiSpadaro}, the definition above follows the approach of~\cite{AmbrosioBV} and~\cite{Reshetnyak}; several other definitions of Sobolev spaces of maps with values in a metric space are available, starting from~\cite{KorevaarSchoen}: their equivalence with the one adopted here was established, partially, in~\cite{Reshetnyak2} and~\cite{Jost} and, in full generality, in~\cite{Chiron}. For $1\leqslant p<+\infty$, we denote by $L^p(\Omega,N)$ the set of the measurable maps $u:\Omega\to N$ such that $d^N(u(\cdot),y_0)\in L^p(\Omega)$ for every $y_0\in N$ and we endow it with the intrinsic distance
$$
d_{L^p}(u,v)=\left(\int_\Omega d^N\big(u(x),v(x)\big)^p\,d\mu_g(x)\right)^{1/p}.
$$
Since $W^{1,p}(\Omega,N)$ is not a vector space when the target $N$ is nonlinear, there is no canonical Banach--space weak topology on it. We therefore adopt the following intrinsic notions of convergence.

Following~\cite{FocardiSpadaro}, for $1<p<+\infty$ we say that a sequence $u_j\in W^{1,p}(\Omega,N)$ converges weakly to $u\in W^{1,p}(\Omega,N)$ if
$$
d_{L^p}(u_j,u)\to0\qquad\text{ and }\qquad\sup_{j\in\NN}\Vert du_j\Vert_{L^p(\Omega)}<+\infty.
$$

For $p=1$, we say that a sequence $u_j\in W^{1,1}(\Omega,N)$ converges weakly to $u\in W^{1,1}(\Omega,N)$ if
$$
d_{L^1}(u_j,u)\to0
$$
and the family $|du_j|$ is equiintegrable in $L^1(\Omega)$. Thus, at $p=1$ our definition is stronger than the one used in~\cite{FocardiSpadaro}, where only uniform boundedness of the differentials in $L^1$ is required. When $N=\R^m$, equiintegrability of the differentials already follows from the usual weak convergence in $W^{1,1}(\Omega,\R^m)$ by the Dunford--Pettis theorem. The equivalence with that convergence at $p=1$ is explained below.

Finally, we say that a sequence $u_j\in W^{1,\infty}(\Omega,N)$ converges weakly$^*$ to $u\in W^{1,\infty}(\Omega,N)$ if
$$
\Vert d^N(u_j,u)\Vert_{L^\infty(\Omega)}\to0\qquad\text{ and }\qquad\sup_{j\in\NN}\Vert du_j\Vert_{L^\infty(\Omega)}<+\infty.
$$
This is again exactly the notion used by Focardi and Spadaro in~\cite{FocardiSpadaro}.

\begin{oss}\label{weakterminology}
The terminology just introduced involves a slight abuse. Since $W^{1,p}(\Omega,N)$ is nonlinear, the expressions ``weak convergence'' and ``weak$^*$ convergence'' do not, in general, refer to the weak or weak$^*$ topology of a Banach space. Throughout the paper, however, for simplicity and unless explicitly stated otherwise, these expressions, as well as the symbols $u_j\rightharpoonup u$ and $u_j\weakkk u$, will always refer to the notions of convergence defined above.

A standard alternative is the extrinsic approach: one chooses an isometric embedding of $N$ into a Euclidean space and uses the ordinary weak or weak$^*$ convergence of the embedded maps. We prefer the intrinsic formulation above, in the same spirit as Focardi and Spadaro~\cite{FocardiSpadaro}, because it does not introduce an ambient Euclidean space and is closer to the theories of Sobolev maps with values in metric spaces recalled above. In particular, the definition of $W^{1,p}(\Omega,N)$ and the convergence notions themselves are largely metric in nature and are therefore better suited to possible extensions to metric targets. On the other hand, the R--quasiconvexity developed in this paper still uses the linear tangent spaces of $M$ and $N$ and the exponential maps; extending the whole theory to general metric spaces would consequently require an appropriate replacement for this tangent structure.
\end{oss}

When $N=\R^m$, the notions of convergence defined above imply the usual weak and weak$^*$ convergences in the corresponding Sobolev spaces, and they are equivalent to them whenever the embedding $W^{1,p}(\Omega,\R^m)\allowbreak\hookrightarrow L^p(\Omega,\R^m)$ is compact, for instance when the bounded open set $\Omega$ has a Lipschitz boundary. For $p=1$, this equivalence holds on every bounded open set: usual weak convergence gives equiintegrability also of the maps, so local compactness yields strong $L^1$ convergence on $\Omega$. For $1<p\leqslant+\infty$, on an arbitrary bounded open set the converse implication holds locally, on every relatively compact smooth subdomain. Consequently, for the nonnegative integrands considered below, requiring lower semicontinuity on every bounded open set with respect to the convergence we adopted is equivalent, in the Euclidean case, to requiring lower semicontinuity with respect to the usual weak (or weak$^*$) convergence. Indeed, if $u_j$ converges to $u$ in the usual weak sense in $W^{1,p}(\Omega,\R^m)$ (or in the usual weak$^*$ sense when $p=+\infty$), choosing an increasing exhaustion $\Omega_k\subseteq\Omega$ by smooth relatively compact open subsets, we have that on each $\Omega_k$, the sequence $u_j$ converges in the sense above, hence
$$
\int_{\Omega_k}f(du)\,d\mu_g\leqslant\liminf_{j\to\infty}\int_{\Omega_k}f(du_j)\,d\mu_g
\leqslant\liminf_{j\to\infty}\int_{\Omega}f(du_j)\,d\mu_g.
$$
Letting $k\to\infty$, the desired inequality on $\Omega$ follows by the monotone convergence theorem. In particular, in the Euclidean case the lower semicontinuity statements below recover exactly the classical Acerbi--Fusco formulation also for $p=1$.

\begin{oss}\label{convFS}
When $1<p\leqslant+\infty$, our notion of convergence coincides exactly with that of Focardi and Spadaro~\cite{FocardiSpadaro}, adapted to our setting where both the domain and the target are Riemannian manifolds. The difference concerns only the case $p=1$, where we also require equiintegrability of the differentials, whereas they require only strong $L^1$ convergence of the maps and uniform boundedness of their differentials in $L^1$. Under this weaker convergence, Theorem~0.6 in~\cite{FocardiSpadaro} states sequential lower semicontinuity of the integral functional for continuous quasiconvex integrands satisfying their growth assumptions, while in Section~\ref{AFR}, we allow Carath\'eodory integrands. The role of equiintegrability is clarified in Remark~\ref{FSpone}, after the proof of Theorem~\ref{AcerbiFuscoRp}: we will give a counterexample showing that this assumption cannot in general be dropped for Carath\'eodory integrands.\\
As for the necessity of the R--quasiconvexity of the integrand, our stronger convergence at $p=1$ makes the assumption of sequential lower semicontinuity of the associated functional weaker, since it is tested along a smaller class of sequences. In both settings, however, weak$^*$ convergence in $W^{1,\infty}$ implies the corresponding convergence in $W^{1,p}$ for every $1\leqslant p<+\infty$. Hence, lower semicontinuity in $W^{1,p}$ implies lower semicontinuity in $W^{1,\infty}$ and the argument proving necessity in $W^{1,\infty}$ also gives the R--quasiconvexity of the integrand for finite $p$.
\end{oss}

\begin{Teorema}\label{QCR->sci}
Let $(M,g)$ and $(N,h)$ be a pair of smooth, connected and complete Riemannian manifolds and let $\mu_g$ be the Riemannian volume measure of $(M,g)$. Let $f:\mathscr{L}(TM,TN)\to \mathbb{R}$ be continuous and R--quasiconvex in the sense of Definition~\ref{RiemannianQC2}. Then, for every bounded open subset $\Omega\subseteq M$, the functional 
$$
u\mapsto F(u,\Omega)=\int_\Omega f(du)\,d\mu_g
$$
is sequentially lower semicontinuous with respect to the weak$^*$ convergence in $W^{1,\infty}(\Omega,N)$, that is,
$$
F(u,\Omega)=\int_{\Omega}f\big(du\big)\,d\mu_g\leqslant\liminf_{j\to\infty }\int_{\Omega}f\big(du_j\big)\,d\mu_g=\liminf_{j\to\infty}F(u_j, \Omega),
$$
for every sequence $u_j\weakkk u$ in $W^{1,\infty}(\Omega,N)$.
\end{Teorema}
\begin{proof}
Let $n\in\NN$ be the dimension of $M$. Let $x_0\in\Omega$ be a point where $u$ is differentiable and where the R--quasiconvexity inequality~\eqref{quasidefR} holds (that is, $\mu_g$--almost~every point of $\Omega$), set $y_0=u(x_0)$ and let us assume that $r>0$ is small enough for $Q_{x_0}^r$ to be contained in a neighborhood $U\subseteq\Omega$ of $x_0$ and for $u(Q_{x_0}^r)$ to be contained in a neighborhood $V=P_{y_0}^s\subseteq N$ of $y_0$ such that the vector bundle $\mathscr{L}(TU,TV)$ can be ``trivialized'' as we discussed in Section~\ref{quasisec}, by choosing orthonormal frames $E_1,\dots,E_n$ in $U$ and $F_1,\dots,F_m$ in $V$, with an associated map $I$ and distance $\delta$ as in formulas~\eqref{mapI} and~\eqref{distL}. Since $u_j\weakkk u$ in $W^{1,\infty}(\Omega,N)$, the functions $u_j$ converge uniformly to $u$ on $Q_{x_0}^r$; hence, by disregarding a finite number of indices, we can assume that $u_j(Q_{x_0}^r)\subseteq V$ for every $j\in\NN$. Moreover, since the maps $u_j$ converge uniformly to $u$, whose image is bounded (as $d^N(u(\cdot),y_0)\in L^\infty(\Omega)$ for every $y_0\in N$), hence relatively compact in the complete manifold $N$ by the Hopf--Rinow theorem, the images $u_j(\Omega)$ are contained in a compact subset of $N$ independent of $j\in\NN$ and the norms $\Vert du_j\Vert_\infty$ are uniformly bounded, all the arguments of the continuous function $f$ appearing below lie in a compact subset of $\mathscr{L}(TM,TN)$, hence $f$ is bounded on it and, by possibly adding a constant to $f$, which changes neither the hypotheses nor the conclusion, we can also assume that $f$ is positive on it. Furthermore, on the compact subset $K$ of $\mathscr{L}(TU,TV)$ containing the arguments of $f$ appearing in the computations on $Q_{x_0}^r$, the function $f$ has a uniform modulus of continuity $\omega$ (with respect to the distance $\delta$), which we can assume to be continuous, bounded and concave, hence subadditive.
We set
$$
v=\exp_{y_0}^{-1}\circ u,\qquad v_j=\exp_{y_0}^{-1}\circ u_j,
$$
and
$$
A[x]=dv[x]=d\exp_{y_0}^{-1}[u(x)]\circ du[x],\qquad A_j[x]=dv_j[x]=d\exp_{y_0}^{-1}[u_j(x)]\circ du_j[x].
$$
We consider a smooth function $\psi:Q_{x_0}^r\to[0,1]$ with compact support and we set, for every $j\in\NN$,
$$
\varphi_j=\exp_{y_0}\big(\psi(v_j-v)\big).
$$
By disregarding again a finite number of indices if necessary, $\varphi_j\in W^{1,\infty}(Q_{x_0}^r,P_{y_0}^s)$ and $\varphi_j=y_0$ outside a compact subset of $Q_{x_0}^r$, hence
\begin{equation}\label{eq590}
d\exp_{y_0}^{-1}[\varphi_j(x)]\circ d\varphi_j[x]=\psi\,d(v_j-v)[x]+d\psi[x]\otimes(v_j-v)=\psi(A_j[x]-A[x])+d\psi[x]\otimes(v_j-v).
\end{equation}
To simplify the notation in the computations below, we define
\begin{equation*}
L_x=d\exp_{x_0}[\exp_{x_0}^{-1}(x)],
\end{equation*}
for every $x\in Q_{x_0}^r$ and we start by applying the hypothesis of R--quasiconvexity for $f$, observing that inequality~\eqref{quasidefR} also holds for maps in $W^{1,\infty}(Q_{x_0}^r,P_{y_0}^s)$ which are equal to $y_0$ outside a compact subset of $Q_{x_0}^r$, by a standard approximation argument in normal coordinates in the target (mollifying $\exp_{y_0}^{-1}\circ\varphi$ and composing
back with $\exp_{y_0}$). Indeed, the function $o(1)$ depends, by hypothesis, in a monotonically nondecreasing way, on the $L^\infty$ norm of the differential of such maps and this approximation can be chosen in such a way as to control the $L^\infty$ norms of the differentials of the approximating maps with the $L^\infty$ norm of the differential of the target map.\\
Hence, setting $\alpha_{x_0}^{y_0}=du[x_0]$ in inequality~\eqref{quasidefR}, for every $j\in\NN$, we have
\begin{align*}
f\big(du[x_0]\big)\leqslant&\,\fint_{Q_{x_0}^r}\!\!\!\!f\big(du[x_0]+d\exp_{y_0}^{-1}[\varphi_j(x)]\circ d\varphi_j[x]\circ\!L_x\big)\,d\mu_g(x)+o_j(1)\\
=&\,\fint_{Q_{x_0}^r}\!\!\!\!f\big(du[x_0]+\psi\,d(v_j-v)[x]\circ\!L_x+d\psi[x]\otimes(v_j-v)\circ\!L_x\big)\,d\mu_g(x)+o_j(1)\\
=&\,\fint_{Q_{x_0}^r}\!\!\!\!f\big(du[x_0]-\psi A[x]\circ\!L_x+\psi A_j[x]\circ\!L_x+d\psi[x]\otimes(v_j-v)\circ\!L_x\big)\,d\mu_g(x)+o_j(1),
\end{align*}
where the function $o_j(1)$ depends monotonically on $\Vert d\varphi_j\Vert_\infty$, as in Definition~\ref{RiemannianQC2}.\\
Then, by the properties of the modulus of continuity $\omega$ of $f$, there holds
\begin{align*}
f\big(du[x_0]-&\,\psi A[x]\circ\!L_x+\psi A_j[x]\circ\!L_x+d\psi[x]\otimes(v_j-v)\circ\!L_x\big)-f\big(A_j[x]\circ\!L_x\big)\\
\leqslant&\,\omega\big(\big|du[x_0]-\psi A[x]\circ\!L_x+(\psi-1)A_j[x]\circ\!L_x+d\psi[x]\otimes(v_j-v)\circ\!L_x\big|_{g_{x_0},h_{y_0}}\big)\\
\leqslant&\,\omega\big(\big|du[x_0]-\psi A[x]\circ\!L_x\big|_{g_{x_0},h_{y_0}}\big)+\omega\big(\big|(\psi-1)A_j[x]\circ\!L_x\big|_{g_{x_0},h_{y_0}}\big)\\
&\,+\omega\big(\big|d\psi[x]\otimes(v_j-v)\circ\!L_x\big|_{g_{x_0},h_{y_0}}\big)\\
\leqslant&\,\omega\big(\big|du[x_0]-\psi A[x]\circ\!L_x\big|_{g_{x_0},h_{y_0}}\big)+\omega\big(C|\psi-1|\big)+\omega\big(\big|d\psi[x]\otimes(v_j-v)\circ\!L_x\big|_{g_{x_0},h_{y_0}}\big),
\end{align*}
where we took into account that the distance $\delta$ on $\mathscr{L}(TU,TV)$, restricted to the fibers, coincides with the one induced by the metrics $g$ and $h$. Thus, we obtain
\begin{align}
f\big(du[x_0]\big)\leqslant&\,\fint_{Q_{x_0}^r}f\big(A_j[x]\circ\!L_x\big)\,d\mu_g(x)+\fint_{Q_{x_0}^r}\omega\big(\big|du[x_0]-\psi A[x]\circ\!L_x\big|_{g_{x_0},h_{y_0}}\big)\,d\mu_g(x)\\
&\,+\fint_{Q_{x_0}^r}\omega\big(C|\psi-1|\big)\,d\mu_g(x)+\fint_{Q_{x_0}^r}\omega\big(\big|d\psi[x]\otimes(v_j-v)\circ\!L_x\big|_{g_{x_0},h_{y_0}}\big)\,d\mu_g(x)+o_j(1).\label{eq556}
\end{align}
We then deal with the term
$$
\fint_{Q_{x_0}^r}f\big(A_j[x]\circ\!L_x\big)\,d\mu_g(x)\leqslant\fint_{Q_{x_0}^r}f\big(du_j[x]\big)\,d\mu_g(x)+\fint_{Q_{x_0}^r}\omega\big(\delta\big(A_j[x]\circ\!L_x,du_j[x]\big)\big)\,d\mu_g(x),
$$
showing that the last integral is bounded by a function $o(1)$ independent of $j\in\NN$, after disregarding a finite number of indices.\\
Indeed, recalling formula~\eqref{distL}, we have
\begin{align}
\omega\big(\delta\big(A_j[x]\circ\!L_x,du_j[x]\big)\big)=&\,\omega\big(d^M(x_0,x)+d^N(y_0,u_j(x))+\big\Vert I(A_j[x]\circ\!L_x)-Idu_j[x]\big\Vert\big)\\
\leqslant&\,\omega\big(d^M(x_0,x)\big)+\omega\big(d^N(y_0,u_j(x))\big)+\omega\big(\big\Vert I(A_j[x]\circ\!L_x)-Idu_j[x]\big\Vert\big)\\
\leqslant&\,\omega(\sqrt n\,r/2)+\omega\big(d^N(y_0,u_j(x))\big)+\omega\big(\big\Vert I(A_j[x]\circ\!L_x)-Idu_j[x]\big\Vert\big).
\end{align}
Let $B(x)$ be the matrix of $L_x:T_{x_0}M\to T_xM$ with respect to the bases given by the frame $E_1,\dots,E_n$ at $x_0$ and $x$, and let $\Lambda(y)$ be the matrix of $d\exp_{y_0}^{-1}[y]:T_yN\to T_{y_0}N$ with respect to the bases given by the frame $F_1,\dots,F_m$ at $y$ and $y_0$. Then
$$
I(A_j[x]\circ\!L_x)=\Lambda(u_j(x))Idu_j[x]B(x),
$$
and hence
\begin{align}
\big\Vert I(A_j[x]\circ\!L_x)-Idu_j[x]\big\Vert\leqslant&\,\big\Vert \Lambda(u_j(x))Idu_j[x](B(x)-\mathrm{Id})\big\Vert+\big\Vert(\Lambda(u_j(x))-\mathrm{Id})Idu_j[x]\big\Vert\\
\leqslant&\,C\Vert du_j\Vert_\infty\big(\Vert B(x)-\mathrm{Id}\Vert+d^N(y_0,u_j(x))\big)=o(1),\label{eq570}
\end{align}
where the last $o(1)$ is independent of $j$, after disregarding a finite number of indices (depending on $r$, which is harmless, since below we let $j\to\infty$ for every fixed $r$): indeed, the norms $\Vert du_j\Vert_\infty$ are uniformly bounded, $B(x)\to\mathrm{Id}$ uniformly as $r\to0$ and $d^N(y_0,u_j(x))\leqslant d^N(u(x_0),u(x))+d^N(u(x),u_j(x))\leqslant C\sqrt n\,r/2+\Vert d^N(u,u_j)\Vert_{L^\infty(Q_{x_0}^r)}$, as $u$ is Lipschitz and $u_j\to u$ uniformly.\\
Hence, inequality~\eqref{eq556} becomes
\begin{align}
f\big(du[x_0]\big)\leqslant&\,\fint_{Q_{x_0}^r}f\big(du_j[x]\big)\,d\mu_g(x)+\fint_{Q_{x_0}^r}\omega\big(\big|du[x_0]-\psi A[x]\circ\!L_x\big|_{g_{x_0},h_{y_0}}\big)\,d\mu_g(x)+o(1)\\
&\,+\fint_{Q_{x_0}^r}\omega\big(C|\psi-1|\big)\,d\mu_g(x)+\fint_{Q_{x_0}^r}\omega\big(\big|d\psi[x]\otimes(v_j-v)\circ\!L_x\big|_{g_{x_0},h_{y_0}}\big)\,d\mu_g(x)+o_j(1)\\
\ &\ \label{eq557}
\end{align}
where the function $o(1)$ is independent of $j\in\NN$, while the functions $o_j(1)$ equal $\eta(r,\Vert d\varphi_j\Vert_\infty)$ for some function $\eta$ going to zero as $r\to0$ and monotone nondecreasing (and we can also clearly assume continuous from the right) in its second argument.\\
Observing now that, by equation~\eqref{eq590} and the smoothness of $\exp_{y_0}$, there holds
\begin{equation}
\Vert d\varphi_j\Vert_\infty\leqslant C\big(\Vert\psi\Vert_\infty\Vert A_j-A\Vert_\infty+\Vert d\psi\Vert_\infty\Vert v_j-v\Vert_\infty\big),
\end{equation}
we have
\begin{align}
o_j(1)=\eta(r,\Vert d\varphi_j\Vert_\infty)&\,\leqslant\eta\big(r,C\Vert\psi\Vert_\infty\Vert A_j-A\Vert_\infty+C\Vert d\psi\Vert_\infty\Vert v_j-v\Vert_\infty\big)\\
&\,\leqslant\eta\big(r,C+C\Vert d\psi\Vert_\infty\Vert v_j-v\Vert_\infty\big),\label{eq595}
\end{align}
for a possibly larger constant $C$, as the norms $\Vert A_j-A\Vert_\infty$ are uniformly bounded. It follows that
$$
\limsup_{j\to\infty}o_j(1)\leqslant\limsup_{j\to\infty}\eta\big(r,C+C\Vert d\psi\Vert_\infty\Vert v_j-v\Vert_\infty\big)=\eta(r,C)=o(1),
$$
with a function $o(1)$ independent of $j\in\NN$, by the properties of the function $\eta$ and the fact that $\Vert v_j-v\Vert_\infty\to0$, as $u_j\to u$ uniformly. Then, passing to the liminf as $j\to\infty$ in formula~\eqref{eq557}, we obtain
\begin{align}
f\big(du[x_0]\big)\leqslant&\,\liminf_{j\to\infty}\fint_{Q_{x_0}^r}f\big(du_j[x]\big)\,d\mu_g(x)+\fint_{Q_{x_0}^r}\omega\big(\big|du[x_0]-\psi A[x]\circ\!L_x\big|_{g_{x_0},h_{y_0}}\big)\,d\mu_g(x)\\
&\,+\fint_{Q_{x_0}^r}\omega\big(C|\psi-1|\big)\,d\mu_g(x)+o(1),\label{eqcar9995}
\end{align}
as the last integral in formula~\eqref{eq557} goes to zero, by the dominated convergence theorem. Now letting $\psi$ go to the characteristic function of $Q_{x_0}^r$, the last integral in this formula vanishes (again by the dominated convergence theorem), hence
$$
f\big(du[x_0]\big)\leqslant\liminf_{j\to\infty}\fint_{Q_{x_0}^r}f\big(du_j[x]\big)\,d\mu_g(x)+\fint_{Q_{x_0}^r}\omega\big(\big|du[x_0]-A[x]\circ\!L_x\big|_{g_{x_0},h_{y_0}}\big)\,d\mu_g(x)+o(1).
$$
Finally, we show that the second integral on the right is bounded by a function $o(1)$. By arguing as above, when we dealt with $A_j$, we have
\begin{align*}
\omega\big(\big|du[x_0]-A[x]\circ\!L_x\big|_{g_{x_0},h_{y_0}}\big)\leqslant&\,\omega\big(\delta\big(du[x_0],du[x]\big)+\delta\big(du[x],A[x]\circ\!L_x\big)\big)\\
\leqslant&\,\omega\big(\delta\big(du[x_0],du[x]\big)\big)+\omega\big(\delta\big(du[x],A[x]\circ\!L_x\big)\big)\\
\leqslant&\,\omega\big(d^M(x_0,x)\big)+\omega\big(d^N(y_0,u(x))\big)+\omega\big(\big\Vert Idu[x_0]-Idu[x]\big\Vert\big)+o(1)\\
\leqslant&\,\omega\big(\big\Vert Idu[x_0]-Idu[x]\big\Vert\big)+o(1),
\end{align*}
hence,
$$
\fint_{Q_{x_0}^r}\omega\big(\big|du[x_0]-A[x]\circ\!L_x\big|_{g_{x_0},h_{y_0}}\big)\,d\mu_g(x)\leqslant\fint_{Q_{x_0}^r}\omega\big(\big\Vert Idu[x_0]-Idu[x]\big\Vert\big)\,d\mu_g(x)+o(1).
$$
If now $x_0$ is a Lebesgue point for the function $Idu$, the integral on the right goes to zero as $r\to0$, by Jensen's inequality, since the function $\omega$ is concave and continuous with $\omega(0)=0$. Thus, we conclude that at $\mu_g$--almost~every point $x_0$ of $\Omega$, there holds
\begin{equation}\label{eq558}
f\big(du[x_0]\big)\leqslant\liminf_{j\to\infty}\fint_{Q_{x_0}^r}f\big(du_j\big)\,d\mu_g+o(1).
\end{equation}
As $\Omega$ is bounded and the functions $f(du_j)$ are uniformly bounded, the following quantity is finite,
\begin{equation*}
\liminf_{j\to\infty}\int_{\Omega}f\big(du_j\big)\,d\mu_g=\ell<+\infty
\end{equation*}
and we can define the finite positive Radon measures $\nu_j$ on $\Omega$, for every $j\in\NN$, given by
\begin{equation*}
d\nu_j=f\big(du_j\big)\,d\mu_g,
\end{equation*}
satisfying the uniform bound $\Vert\nu_j\Vert=\int_\Omega f\big(du_j\big)\,d\mu_g\leqslant C$. Hence, up to a subsequence, we may assume that $\int_\Omega f(du_j)\,d\mu_g\to\ell$ and, by the Banach--Alaoglu--Bourbaki theorem, that $\nu_j$ weakly$^*$ converges to some limit Radon measure $\nu$.\\
For every open subset $G$ of $\Omega$, by the properties of the weak$^*$ convergence, we have
\begin{equation*}
\nu(G)\leqslant\liminf_{j\to\infty}\nu_j(G)=\liminf_{j\to\infty}\int_G f\big(du_j\big)\,d\mu_g\leqslant C\mu_g(G),
\end{equation*}
where $C$ is an upper bound for the functions $f(du_j)$, hence, by the outer regularity of the measure $\mu_g$, there holds $\nu(B)\leqslant C\mu_g(B)$ for every Borel set $B\subseteq\Omega$, in particular $\nu\ll\mu_g$. Thus, we may apply the Radon--Nikodym theorem, obtaining a function $\theta:\Omega\to[0,+\infty)$ such that $d\nu=\theta\,d\mu_g$ and, by the differentiation theorem for Radon measures (in normal coordinates around $x$, the closed cubes $\overline Q_x^r$ form a ``regular'' family),
\begin{equation*}
\theta(x)=\lim_{r\to0}\frac{\nu\big(\overline Q_x^r\big)}{\mu_g\big(\overline Q_x^r\big)},
\end{equation*}
for $\mu_g$--almost~every $x\in\Omega$.\\
Then, at the points $x_0\in\Omega$ where inequality~\eqref{eq558} is satisfied and the above limit holds, we have
\begin{align*}
f\big(du[x_0]\big)\leqslant&\,\liminf_{r\to0}\bigg(\liminf_{j\to\infty}\fint_{Q_{x_0}^r}f\big(du_j\big)\,d\mu_g+o(1)\bigg)=\liminf_{r\to0}\liminf_{j\to\infty}\frac{\nu_j\big(Q^r_{x_0}\big)}{\mu_g\big(Q^r_{x_0}\big)}\\
=&\,\liminf_{r\to0}\liminf_{j\to\infty}\frac{\nu_j\big(\overline Q^r_{x_0}\big)}{\mu_g\big(\overline Q^r_{x_0}\big)}\leqslant\liminf_{r\to0}\frac{\nu\big(\overline Q^r_{x_0}\big)}{\mu_g\big(\overline Q^r_{x_0}\big)}=\theta(x_0),
\end{align*}
as $\limsup_{j\to\infty}\nu_j\big(\overline Q^r_{x_0}\big)\leqslant\nu\big(\overline Q^r_{x_0}\big)$, since the sets $\overline Q^r_{x_0}$ are compact. Hence, $f\big(du[x]\big)\leqslant\theta(x)$ $\mu_g$--almost~everywhere in $\Omega$, thus
\begin{align*}
F(u,\Omega)&=\int_\Omega f\big(du\big)\,d\mu_g\leqslant\int_\Omega\theta\,d\mu_g=\nu(\Omega)\leqslant\liminf_{j\to\infty}\nu_j(\Omega)\\
&=\liminf_{j\to\infty}\int_\Omega f\big(du_j\big)\,d\mu_g=\liminf_{j\to\infty}F(u_j,\Omega)
\end{align*}
and the proof is complete.
\end{proof}

The following theorem is the converse of the previous one.

\begin{Teorema}\label{sci->QCR}
Let $(M,g)$ and $(N,h)$ be a pair of smooth, connected and complete Riemannian manifolds and let $\mu_g$ be the Riemannian volume measure of $(M,g)$. Let $f:\mathscr{L}(TM,TN)\to \mathbb{R}$ be a continuous function. If for every bounded open subset $\Omega\subseteq M$, the functional
\begin{equation*}
F(u,\Omega)=\int_{\Omega}f\big(du\big)\,d\mu_g
\end{equation*}
is sequentially lower semicontinuous with respect to the weak$^*$ convergence in $W^{1,\infty}(\Omega,N)$, then the function $f$ is R--quasiconvex in the sense of Definition~\ref{RiemannianQC2}. 
\end{Teorema}
\begin{proof}
Let $n\in\NN$ be the dimension of $M$. We adopt the same setting and notation as in the proof of Theorem~\ref{QCR->sci}. In particular, having fixed $x_0\in M$, $y_0\in N$ and $\alpha_{x_0}^{y_0}\in\mathscr L(T_{x_0}M,T_{y_0}N)$, we work in neighborhoods $U\subseteq M$ of $x_0$ and $V=P_{y_0}^s\subseteq N$ of $y_0$ such that the vector bundle $\mathscr L(TU,TV)$ can be ``trivialized'' 
and $Q_{x_0}^r\subseteq U$.\\
We consider a smooth function $u:U\to N$ such that $u(x_0)=y_0$ and $du[x_0]=\alpha_{x_0}^{y_0}$ (which clearly exists, after possibly reducing $U$), with $u(U)\subseteq P_{y_0}^{s/2}$, and set $v=\exp_{y_0}^{-1}\circ u$. Let $\varphi\in C^{\infty}(Q_{x_0}^r,P_{y_0}^s)$ be equal to $y_0$ outside a compact subset of $Q_{x_0}^r$ and set $\widetilde\varphi=\exp_{y_0}^{-1}\circ\varphi$. Then the function $\widetilde\varphi\circ\exp_{x_0}$ is defined on the open cube $\widetilde Q_{x_0}^r=r(-1/2,1/2)^n\subseteq T_{x_0}M$, smooth and with compact support. We define the extension by periodicity $\psi:T_{x_0}M\to T_{y_0}N$ of $\widetilde\varphi\circ\exp_{x_0}$ to the whole of $T_{x_0}M$ and we set
$$
\varphi_k(x)=\frac{1}{k}\psi(k\exp_{x_0}^{-1}(x))
$$
for every $x\in Q_{x_0}^r$ and $k\in\NN$. Then, clearly $\varphi_k\to0$ uniformly and the differentials $d\varphi_k$ are uniformly bounded, since
$$
d\varphi_k[x]=d\psi\big[k\exp_{x_0}^{-1}(x)\big]\circ d\exp_{x_0}^{-1}[x].
$$
Hence, for $k\in\NN$ large enough, we can define
$$u_k(x)=\exp_{y_0}\big(v(x)+\varphi_k(x)\big),$$
and $u_k\weakkk u$ in $W^{1,\infty}(Q_{x_0}^r,N)$, since $u_k\to u$ uniformly and the differentials $du_k$ are uniformly bounded. Thus, by the hypothesis of lower semicontinuity of the functional $F$, there holds
\begin{equation}\label{eq571N}
\int_{Q_{x_0}^r}f\big(du[x]\big)\,d\mu_g(x)\leqslant\liminf_{k\to\infty}\int_{Q_{x_0}^r}f\big(du_k[x]\big)\,d\mu_g(x).
\end{equation}
For $x\in Q_{x_0}^r$, set $z=\exp_{x_0}^{-1}(x)$. Since
$$
du_k[x]=d\exp_{y_0}\big[v(x)+\varphi_k(x)\big]\circ\Big(dv[x]+d\psi[kz]\circ d\exp_{x_0}^{-1}[x]\Big),
$$
arguing as in the proof of Theorem~\ref{QCR->sci} and using the smoothness of $u$, of the exponential maps and of the local orthonormal frames, we have, uniformly for $x\in Q_{x_0}^r$,
\begin{equation}\label{eq580a}
\delta\Big(du_k[x],\alpha_{x_0}^{y_0}+d\psi[kz]\Big)\leqslant Cr\big(1+\Vert d\psi\Vert_\infty\big)+\frac{C}{k}\Vert\psi\Vert_\infty\big(1+\Vert d\psi\Vert_\infty\big),
\end{equation}
for $r>0$ small enough and $k\in\NN$ large enough. Indeed, representing the linear maps in the chosen local orthonormal frames, $dv[x]=\alpha_{x_0}^{y_0}+O(r)$, $d\exp_{x_0}^{-1}[x]=\mathrm{Id}_{T_{x_0}M}+O(r)$ and $d\exp_{y_0}[v(x)+\varphi_k(x)]=\mathrm{Id}_{T_{y_0}N}+O(r+k^{-1}\Vert\psi\Vert_\infty)$, while $d^M(x,x_0)=O(r)$ and $d^N(u_k(x),y_0)=O(r+k^{-1}\Vert\psi\Vert_\infty)$. Therefore, considering a suitable continuous, bounded and concave modulus of continuity $\omega$ of $f$, formula~\eqref{eq580a} gives
\begin{align}
\big|f\big(du_k[x]\big)-f\big(\alpha_{x_0}^{y_0}+d\psi[k\exp_{x_0}^{-1}(x)]\big)\big|\leqslant&\,\omega\Big(Cr(1+\Vert d\psi\Vert_\infty)\Big)+\omega\bigg(\frac{C}{k}\Vert\psi\Vert_\infty(1+\Vert d\psi\Vert_\infty)\bigg)\\
\leqslant&\,o(1)(1+\Vert d\psi\Vert_\infty)+o_k(1),\label{eq580}
\end{align}
as $r\to0$ and $k\to\infty$, where $o_k(1)\to0$ as $k\to\infty$ for fixed $r$. Hence, from~\eqref{eq571N},
\begin{equation}\label{eq572}
\int_{Q_{x_0}^r}f\big(du[x]\big)\,d\mu_g(x)\leqslant\liminf_{k\to\infty}\int_{Q_{x_0}^r}f\big(\alpha_{x_0}^{y_0}+d\psi\big[k\exp_{x_0}^{-1}(x)\big]\big)\,d\mu_g(x)+o(1)\mu_g(Q_{x_0}^r).
\end{equation}
Changing variables as $z=\exp_{x_0}^{-1}(x)$, so $z\in r(-1/2,1/2)^n\subseteq T_{x_0}M$ and denoting by $dz$ the Lebesgue measure on $T_{x_0}M$ relative to the metric tensor $g_{x_0}$, we have
\begin{align}
\int_{Q_{x_0}^r}f\big(\alpha_{x_0}^{y_0}+d\psi\big[k\exp_{x_0}^{-1}(x)\big]\big)\,d\mu_g(x)=&\,\int_{r(-1/2,1/2)^n}f\big(\alpha_{x_0}^{y_0}+d\psi[kz]\big)\,J_{\exp_{x_0}}(z)\,dz\\
\leqslant&\,\int_{r(-1/2,1/2)^n}f\big(\alpha_{x_0}^{y_0}+d\psi[kz]\big)\,dz+Cr^{n+1},\label{eq573}
\end{align}
for some constant $C$ independent of $k\in\NN$, since the Jacobian of the exponential map $\exp_{x_0}$ at the point $z\in T_{x_0}M$, namely
$$
J_{\exp_{x_0}}(z)=\Big|\det_{ij}\big(g_{\exp_{x_0}(z)}(d\exp_{x_0}[z](E_i),E_j)\big)\Big|,
$$
satisfies
$$
J_{\exp_{x_0}}(z)=J_{\exp_{x_0}^{-1}}\big(\exp_{x_0}(z)\big)^{-1}.
$$
The function $J_{\exp_{x_0}}$ is smooth and goes uniformly to $1$ as $r\to0$, because $|z|\leqslant r\sqrt n/2$ (we recall that $d\exp_{x_0}[0]$ is the identity of $T_{x_0}M$) and $f\big(\alpha_{x_0}^{y_0}+d\psi[kz]\big)$ is bounded. Then, changing the variables again as $w=kz$ in the last integral, we get
$$
\int_{Q_{x_0}^r}f\big(\alpha_{x_0}^{y_0}+d\psi\big[k\exp_{x_0}^{-1}(x)\big]\big)\,d\mu_g(x)\leqslant\frac{1}{k^n}\int_{kr(-1/2,1/2)^n}f\big(\alpha_{x_0}^{y_0}+d\psi[w]\big)\,dw+Cr^{n+1}.
$$
Now, by the periodicity of $\psi$, there holds
$$
\frac{1}{k^n}\int_{kr(-1/2,1/2)^n}f\big(\alpha_{x_0}^{y_0}+d\psi[w]\big)\,dw=\int_{r(-1/2,1/2)^n}f\big(\alpha_{x_0}^{y_0}+d\widetilde\varphi[\exp_{x_0}(w)]\circ d\exp_{x_0}[w]\big)\,dw
$$
and if we change (back) variables as $x=\exp_{x_0}(w)$, we obtain
\begin{align*}
\int_{r(-1/2,1/2)^n}&\,f\big(\alpha_{x_0}^{y_0}+d\widetilde\varphi[\exp_{x_0}(w)]\circ d\exp_{x_0}[w]\big)\,dw\\
=&\,\int_{Q_{x_0}^r}f\big(\alpha_{x_0}^{y_0}+d\exp_{y_0}^{-1}[\varphi(x)]\circ d\varphi[x]\circ d\exp_{x_0}[\exp_{x_0}^{-1}(x)]\big)\,J_{\exp_{x_0}^{-1}}(x)\,d\mu_g(x)\\
\leqslant&\,\int_{Q_{x_0}^r}f\big(\alpha_{x_0}^{y_0}+d\exp_{y_0}^{-1}[\varphi(x)]\circ d\varphi[x]\circ d\exp_{x_0}[\exp_{x_0}^{-1}(x)]\big)\,d\mu_g(x)+Cr\mu_g(Q_{x_0}^r),
\end{align*}
arguing as above. Hence, by this inequality, whose right--hand side does not depend on $k\in\NN$, and by formulas~\eqref{eq572} and~\eqref{eq573}, we conclude
\begin{align*}
\int_{Q_{x_0}^r}f\big(du[x]\big)\,d\mu_g(x)\leqslant&\,\int_{Q_{x_0}^r}f\big(\alpha_{x_0}^{y_0}+d\exp_{y_0}^{-1}[\varphi(x)]\circ d\varphi[x]\circ d\exp_{x_0}[\exp_{x_0}^{-1}(x)]\big)\,d\mu_g(x)\\
&+o(1)\mu_g(Q_{x_0}^r),
\end{align*}
as $\mu_g(Q_{x_0}^r)\approx r^n$.\\
Since, by the continuity of $du$, there holds
\begin{equation*}
\int_{Q_{x_0}^r}f\big(\alpha_{x_0}^{y_0}\big)\,d\mu_g(x)\leqslant\int_{Q_{x_0}^r}f\big(du[x]\big)\,d\mu_g(x)+o(1)\mu_g(Q_{x_0}^r)
\end{equation*}
and $du[x_0]=\alpha_{x_0}^{y_0}$, we finally have
\begin{align*}
f(\alpha_{x_0}^{y_0})\leqslant&\,\frac{1}{\mu_g(Q_{x_0}^r)}\int_{Q_{x_0}^r}f\big(\alpha_{x_0}^{y_0}+d\exp_{y_0}^{-1}[\varphi(x)]\circ d\varphi[x]\circ d\exp_{x_0}[\exp_{x_0}^{-1}(x)]\big)\,d\mu_g(x)+o(1)
\end{align*}
for every $x_0\in M$, $y_0\in N$ and $\alpha_{x_0}^{y_0}\in\mathscr L(T_{x_0}M,T_{y_0}N)$, with $o(1)$ going to zero as $r\to0$, depending only on the $L^\infty$ norm of $d\varphi$. Moreover, $\Vert d\psi\Vert_\infty\leqslant C\Vert d\varphi\Vert_\infty$, with a constant depending only on the fixed local charts, and by tracing back how we obtained such a function $o(1)$, it is clear that it can be chosen in a way that such a dependence is monotonically nondecreasing (see in particular estimate~\eqref{eq580}). Hence, the function $f$ is R--quasiconvex according to Definition~\ref{RiemannianQC2}.
\end{proof}

Putting together the two theorems, we have the following one, which generalizes (in our case of $f$ continuous) the results of Acerbi and Fusco for $p=+\infty$. 

\begin{Teorema}\label{QCR=sci}
Let $(M,g)$ and $(N,h)$ be a pair of smooth, connected and complete Riemannian manifolds and let $\mu_g$ be the Riemannian volume measure of $(M,g)$. Let $f:\mathscr{L}(TM,TN)\to \mathbb{R}$ be continuous. 
Then, $f$ is R--quasiconvex in the sense of Definition~\ref{RiemannianQC2} if and only if for every bounded open subset $\Omega\subseteq M$, the functional 
$$
u\mapsto F(u,\Omega)=\int_\Omega f(du)\,d\mu_g
$$
is sequentially lower semicontinuous with respect to the weak$^*$ convergence in $W^{1,\infty}(\Omega,N)$.
\end{Teorema}

\begin{oss}
We notice that, if $f$ is continuous and R--quasiconvex, the quasiconvexity condition of Definition~\ref{RiemannianQC2} actually holds for every $x_0\in M$, not only for $\mu_g$--almost~every $x_0\in M$, as indeed the proof of Theorem~\ref{sci->QCR} shows. This can also be seen directly: the set of the points $x_0\in M$ such that $f_{x_0}^{y_0}$ is quasiconvex for every $y_0\in N$ is closed, since the quasiconvexity inequality (the third condition in Proposition~\ref{RiemannianQC-equiv}, written in a local ``trivialization'' of the bundle) passes to the limit as $x_0$ varies, by the uniform continuity of $f$ on the compact subsets of $\mathscr{L}(TM,TN)$; having full measure, it is also dense, hence it coincides with $M$. The distinction between ``every'' and ``almost~every'' $x_0$ becomes relevant only for Carath\'eodory integrands, which are treated in the next section.
\end{oss}

\section{The Riemannian analogues of the Acerbi--Fusco theorems}\label{AFR}

We now extend the results of the previous section in the two directions appearing in the theorems of Acerbi and Fusco recalled in the Introduction: first, we allow a merely measurable dependence of the integrand on the point $x\in M$ and, second, we consider the spaces $W^{1,p}(\Omega,N)$ for $1\leqslant p<+\infty$. The natural class of integrands is the Riemannian analogue of the Carath\'eodory functions in Theorems~\ref{Fusco1} and~\ref{Teocap1pfinito}. We first make this notion precise.

\medskip

Let $U\subseteq M$ and $V\subseteq N$ be open sets endowed with local orthonormal frames as in Section~\ref{quasisec}, and let
$$
\Phi_{U,V}:\mathscr L(TU,TV)\to U\times V\times\R^{m\times n},\qquad \Phi_{U,V}(\alpha_x^y)=(x,y,I\alpha_x^y)
$$
be the corresponding local ``trivialization'' of the bundle $\mathscr L(TM,TN)$. Moreover, we denote by $|\alpha_x^y|$ the norm of $\alpha_x^y\in\mathscr L(T_xM,T_yN)$ induced by the metrics $g_x$ and $h_y$, that is, the norm associated with the distance $\delta$ restricted to the fibers.

\begin{definizione}\label{CaratheodoryR}
We say that a function $f:\mathscr L(TM,TN)\to\R$ is a {\em Carath\'eodory integrand} if, for every pair of open sets $U\subseteq M$, $V\subseteq N$ as above, the function
$$
f_{U,V}(x,y,A)=f\big(\Phi_{U,V}^{-1}(x,y,A)\big)
$$
is measurable with respect to $x\in U$, for every $(y,A)\in V\times\R^{m\times n}$, and continuous with respect to $(y,A)\in V\times\R^{m\times n}$, for $\mu_g$--almost~every $x\in U$.
\end{definizione}

This definition is independent of the choice of the local orthonormal frames, since the transition maps between two local ``trivializations'' are smooth. In particular, covering $M$ and $N$ with countably many such open sets $U$ and $V$, for a Carath\'eodory integrand the restricted maps $f_{x_0}^{y_0}$ are continuous for $\mu_g$--almost~every $x_0\in M$ and every $y_0\in N$, hence Definition~\ref{RiemannianQC2} applies without any modification.

\medskip

We use the intrinsic definition of the Sobolev spaces $W^{1,p}(\Omega,N)$ introduced in the previous section and, according to the convention in Remark~\ref{weakterminology}, the terms weak and weak$^*$ convergence always refer to the notions fixed there. For $1<p<+\infty$ and for $p=+\infty$ they coincide with the notions of Focardi and Spadaro, whereas for $p=1$ we require in addition the equiintegrability of the differentials.

\medskip

Our first result is the full analogue of Theorem~\ref{Fusco1} in the case $p=+\infty$.

\begin{Teorema}[{Riemannian Acerbi--Fusco -- Case $p=+\infty$}]\label{AcerbiFuscoRinf}
Let $(M,g)$ and $(N,h)$ be a pair of smooth, connected and complete Riemannian manifolds and let $\mu_g$ be the Riemannian volume measure of $(M,g)$. Let $f:\mathscr L(TM,TN)\to\R$ be a Carath\'eodory integrand in the sense of Definition~\ref{CaratheodoryR}, satisfying the following growth condition: there exists a nonnegative function $a\in L^1_{\mathrm{loc}}(M)$ such that, for every compact subset $K\subseteq N$ and every $R>0$, there exists a constant $C_{K,R}\geqslant0$ for which
\begin{equation}\label{growthAFinf}
0\leqslant f(\alpha_x^y)\leqslant a(x)+C_{K,R},
\end{equation}
for $\mu_g$--almost~every $x\in M$, every $y\in K$ and every $\alpha_x^y\in\mathscr L(T_xM,T_yN)$ with $|\alpha_x^y|\leqslant R$.\\
Then, $f$ is R--quasiconvex in the sense of Definition~\ref{RiemannianQC2} if and only if for every bounded open subset $\Omega\subseteq M$, the functional
$$
u\mapsto F(u,\Omega)=\int_\Omega f(du)\,d\mu_g
$$
is sequentially lower semicontinuous with respect to the weak$^*$ convergence in $W^{1,\infty}(\Omega,N)$ (in the sense of Section~\ref{RDefSec}).
\end{Teorema}

The growth assumption above is simply the intrinsic version of the condition
$$
0\leqslant f(x,s,\xi)\leqslant a(x)+b(s,\xi)
$$
in Theorem~\ref{Fusco1}, with $b$ locally bounded. Notice that Theorem~\ref{QCR=sci} already gives the same equivalence for arbitrary continuous integrands, without any growth or sign assumption. Thus, in the case $p=+\infty$, the genuinely new point is the merely measurable dependence on $x$ and the growth condition above provides the required locally integrable upper bound in the Carath\'eodory setting. It is then natural to combine the geometric argument of the previous section with the approximation arguments used by Acerbi and Fusco~\cite{AcerbiFusco} (see also Marcellini~\cite{marcellini}).

\begin{oss}\label{crescitacontinua}
The growth condition in Theorem~\ref{AcerbiFuscoRinf} is automatically satisfied by every nonnegative continuous integrand. Indeed, let $M_j\subseteq M$ and $K_j\subseteq N$ be increasing exhaustions by compact subsets, chosen so that $M_j\subseteq\mathop{\rm int}M_{j+1}$ and $K_j\subseteq\mathop{\rm int}K_{j+1}$, and set
$$
T_j(x)=\max\big\{f(\alpha_x^y)\,\,|\,\,y\in K_j\text{ and }|\alpha_x^y|\leqslant j\big\}\qquad\text{ and }\qquad C_j
=\max_{x\in M_j}T_j(x).
$$
For every $j\in\NN$ and every $x\in M$, the set over which the maximum defining $T_j(x)$ is taken is compact, hence $T_j(x)<+\infty$. Moreover, $T_j$ is continuous by the maximum theorem applied in local ``trivializations'' of $\mathscr L(TM,TN)$. Since $M_j\times K_j$ is compact, the corresponding subset with $|\alpha|\leqslant j$ is compact in the bundle, and therefore $C_j=\max_{x\in M_j}T_j(x)<+\infty$. Then, the function $a=\sup_{j\in\NN}(T_j-C_j)^+$ is nonnegative, lower semicontinuous and bounded on every $M_i$, since $T_j\leqslant C_j$ on $M_i$ for every $j\geqslant i$, while $T_j\leqslant T_i\leqslant C_i$ on $M_i$ for every $j<i$. Hence $a\in L^1_{\mathrm{loc}}(M)$ and
$$
f(\alpha_x^y)\leqslant a(x)+C_j
$$
for every $x\in M$, $y\in K_j$ and $|\alpha_x^y|\leqslant j$, which is the growth condition, since every compact subset $K\subseteq N$ is contained in $K_j$ and every $R>0$ satisfies $R\leqslant j$, for $j$ large enough. Therefore, as far as nonnegative integrands are concerned, Theorem~\ref{QCR=sci} is a special case of Theorem~\ref{AcerbiFuscoRinf} and the extra generality of the former lies only in the absence of any sign condition on $f$. By contrast, the growth condition in Theorem~\ref{AcerbiFuscoRp} is a genuine restriction also for continuous integrands, exactly as in Theorem~\ref{Teocap1pfinito}.
\end{oss}

We next consider the case $1\leqslant p<+\infty$. Fixing a point $y_*\in N$, the intrinsic counterpart of the term $|s|^p$ appearing in Theorem~\ref{Teocap1pfinito} is $d^N(y,y_*)^p$.

\begin{Teorema}[Riemannian Acerbi--Fusco -- Case $1\leqslant p<+\infty$]\label{AcerbiFuscoRp}
Let $(M,g)$ and $(N,h)$ be a pair of smooth, connected and complete Riemannian manifolds and let $\mu_g$ be the Riemannian volume measure of $(M,g)$. Let $1\leqslant p<+\infty$, fix $y_*\in N$ and let $f:\mathscr L(TM,TN)\to\R$ be a Carath\'eodory integrand in the sense of Definition~\ref{CaratheodoryR}, satisfying
\begin{equation}\label{growthAFp}
0\leqslant f(\alpha_x^y)\leqslant a(x)+C\big(d^N(y,y_*)^p+|\alpha_x^y|^p\big)
\end{equation}
for $\mu_g$--almost~every $x\in M$, every $y\in N$ and every $\alpha_x^y\in\mathscr L(T_xM,T_yN)$, where $a\in L^1_{\mathrm{loc}}(M)$ is nonnegative and $C\geqslant0$ is a constant.\\
Then, $f$ is R--quasiconvex in the sense of Definition~\ref{RiemannianQC2} if and only if for every bounded open subset $\Omega\subseteq M$, the functional
$$
u\mapsto F(u,\Omega)=\int_\Omega f(du)\,d\mu_g
$$
is sequentially lower semicontinuous with respect to the weak convergence in $W^{1,p}(\Omega,N)$ (in the sense of Section~\ref{RDefSec}).
\end{Teorema}

The growth assumption in the theorem above is the intrinsic counterpart of
$$
0\leqslant f(x,s,\xi)\leqslant a(x)+C\big(|s|^p+|\xi|^p\big)
$$
in Theorem~\ref{Teocap1pfinito}. The choice of the reference point $y_*$ is immaterial, since, by the triangle inequality, changing it only amounts to changing the constant $C$ and adding a constant term, which can be absorbed into the locally integrable function $a$. In particular, no compactness assumption on the target manifold $N$ is needed. We stress that the intrinsic definition of $W^{1,p}(\Omega,N)$ and the notion of convergence considered on this space are two separate issues. For $p=1$, the weak convergence appearing in the theorem is the one defined in Section~\ref{RDefSec}: it requires equiintegrability of the differentials and is
stronger than the convergence used by Focardi and Spadaro. When $N=\R^m$, it coincides with the usual weak convergence in $W^{1,1}(\Omega,\R^m)$ on every bounded open set, as explained in Section~\ref{RDefSec}.

\medskip

When both the domain and the target are Euclidean and $y_*=0$, Theorems~\ref{AcerbiFuscoRinf} and~\ref{AcerbiFuscoRp} reduce, at the level of the lower semicontinuity statements, to Theorems~\ref{Fusco1} and~\ref{Teocap1pfinito}, respectively, taking into account the discussion on the notions of convergence in Section~\ref{RDefSec}. We stress that no generality is lost in this reduction: as shown there, requiring lower semicontinuity with respect to the notions of convergence adopted here for every bounded open subset is equivalent, in the Euclidean case, to requiring it with respect to the usual weak and weak$^*$ convergences. On the other hand, if only the target is Euclidean, these statements extend the results obtained in~\cite{CLM1}, while in the complementary case of a Euclidean domain and a manifold target the sufficiency of the corresponding notion of quasiconvexity for the lower semicontinuity in $W^{1,p}$, for $1\leqslant p\leqslant+\infty$, with respect to the convergence considered there, is stated in Theorem~0.6 of Focardi and Spadaro~\cite{FocardiSpadaro} for continuous integrands. For $p=1$, their result does not require equiintegrability of the differentials, whereas Theorem~\ref{AcerbiFuscoRp} requires this condition but allows Carath\'eodory integrands (see Remark~\ref{convFS}).

\medskip

We prove below the necessity direction in the case $p=+\infty$, by extending to Carath\'eodory integrands the localized oscillation argument used in Theorem~\ref{sci->QCR}, by means of the classical Scorza~Dragoni theorem~\cite{ScorzaDragoni}. The necessity for $1\leqslant p<+\infty$ is then an immediate consequence of the case $p=+\infty$. The main new point is the sufficiency for Carath\'eodory integrands and, for finite $p$, the control of the $p$--growth terms. The proof combines the geometric localization developed in Theorem~\ref{QCR->sci} with the measure--theoretic and approximation arguments of the classical Acerbi--Fusco theory~\cite{AcerbiFusco} and its later variants~\cite{marcellini}. This separates the analytic issues as much as possible from the genuinely Riemannian ones treated in Section~\ref{RDefSec}.

\begin{proof}[Proof of the necessity -- The case $p=+\infty$]\ \\
We only need to modify the proof of Theorem~\ref{sci->QCR} in order to deal with the merely measurable dependence of $f$ on the point $x\in M$. We use the classical Scorza~Dragoni argument~\cite{ScorzaDragoni}, as in the Euclidean case of Acerbi and Fusco~\cite{AcerbiFusco}.\\
Consider a local ``trivialization'' of $\mathscr L(TM,TN)$ as in Section~\ref{quasisec} and write $f_{U,V}(x,y,A)$ for the corresponding expression of $f$. By the Scorza~Dragoni theorem, using a countable atlas of $M$ and $N$, countable exhaustions by compact subsets of the variables $(y,A)$ and compact subsets of the coordinate neighborhoods whose complements have arbitrarily small measure, and then by the Lebesgue density theorem, there exists a set $M_0\subseteq M$ of full $\mu_g$--measure with the following property: for every $x_0\in M_0\cap U$ and every compact subset $K\subseteq V\times\R^{m\times n}$ there exists a compact set $E\subseteq U$ containing $x_0$, of which $x_0$ is a density point, such that $f_{U,V}$ is continuous on $E\times K$. We can also assume that every $x_0\in M_0$ is a Lebesgue point of the function $a$ appearing in the growth assumption, that the growth condition holds at $x_0$ and that $f_{U,V}(x_0,\cdot,\cdot)$ is continuous.\\
We fix $x_0\in M_0$, $y_0\in N$ and $\alpha_{x_0}^{y_0}\in\mathscr L(T_{x_0}M,T_{y_0}N)$, and we follow exactly the proof of Theorem~\ref{sci->QCR}. In particular, we choose $u$, $v$, $\psi$, $\varphi_k$ and $u_k$ as there, so that $u_k\weakkk u$ in $W^{1,\infty}(Q_{x_0}^r,N)$ and
\begin{equation}\label{eqCarNec1}
\int_{Q_{x_0}^r}f\big(du[x]\big)\,d\mu_g(x)\leqslant\liminf_{k\to\infty}\int_{Q_{x_0}^r}f\big(du_k[x]\big)\,d\mu_g(x).
\end{equation}
Moreover, setting $z=\exp_{x_0}^{-1}(x)$, estimate~\eqref{eq580a} gives
\begin{equation}\label{eqCarNec2}
\delta\big(du_k[x],\alpha_{x_0}^{y_0}+d\psi[kz]\big)\leqslant Cr\big(1+\Vert d\psi\Vert_\infty\big)+\frac{C}{k}\Vert\psi\Vert_\infty\big(1+\Vert d\psi\Vert_\infty\big).
\end{equation}
All the arguments of $f$ appearing in this estimate belong, for $r$ small enough and $k$ large enough, to a fixed compact subset $K$ of $V\times\R^{m\times n}$, in the corresponding local ``trivialization''. Let $E$ be a compact set as above, associated with such a set $K$, and set
$$
B_r=Q_{x_0}^r\setminus E.
$$
Since $f_{U,V}$ is uniformly continuous on $E\times K$, from inequality~\eqref{eqCarNec2} we obtain, for $x\in Q_{x_0}^r\cap E$,
\begin{align*}
\big|f\big(du_k[x]\big)-f\big(\alpha_{x_0}^{y_0}+d\psi[k\exp_{x_0}^{-1}(x)]\big)\big|
\leqslant&\,\omega_E\big(Cr(1+\Vert d\psi\Vert_\infty)\big)+\omega_E\Big(\frac{C}{k}\Vert\psi\Vert_\infty(1+\Vert d\psi\Vert_\infty)\Big),
\end{align*}
where $\omega_E$ is a modulus of continuity depending on the compact set $E\times K$. On the other hand, by the growth assumption, for some constant $C$ independent of $r$ and $k$,
\begin{align*}
\int_{B_r}\big|f\big(du_k[x]\big)-f\big(\alpha_{x_0}^{y_0}+d\psi[k\exp_{x_0}^{-1}(x)]\big)\big|\,d\mu_g(x)
\leqslant&\,\int_{B_r}a(x)\,d\mu_g(x)+\big(a(x_0)+C\big)\mu_g(B_r).
\end{align*}
Since $x_0$ is a density point of $E$ and a Lebesgue point of $a$, there holds
$$
\frac{\mu_g(B_r)}{\mu_g(Q_{x_0}^r)}\to0
$$
and
\begin{align*}
\frac{1}{\mu_g(Q_{x_0}^r)}\int_{B_r}a(x)\,d\mu_g(x)
\leqslant&\,\fint_{Q_{x_0}^r}|a(x)-a(x_0)|\,d\mu_g(x)+a(x_0)\frac{\mu_g(B_r)}{\mu_g(Q_{x_0}^r)}\to0
\end{align*}
as $r\to0$. Hence, letting first $k\to\infty$ for every fixed $r$ in inequality~\eqref{eqCarNec1}, we obtain exactly the analogue of formula~\eqref{eq572},
\begin{equation}\label{eqCarNec3}
\int_{Q_{x_0}^r}f\big(du[x]\big)\,d\mu_g(x)\leqslant\liminf_{k\to\infty}\int_{Q_{x_0}^r}f\big(\alpha_{x_0}^{y_0}+d\psi[k\exp_{x_0}^{-1}(x)]\big)\,d\mu_g(x)+o(1)\mu_g(Q_{x_0}^r).
\end{equation}
The same argument, without the oscillating term, gives
\begin{equation}\label{eqCarNec4}
\fint_{Q_{x_0}^r}f\big(du[x]\big)\,d\mu_g(x)=f\big(\alpha_{x_0}^{y_0}\big)+o(1),
\end{equation}
since $u(x_0)=y_0$ and $du[x_0]=\alpha_{x_0}^{y_0}$.\\
Starting from estimate~\eqref{eqCarNec3}, we can now repeat verbatim the argument from formula~\eqref{eq572} to the end of the proof of Theorem~\ref{sci->QCR}: all the arguments of $f$ appearing there belong to the fixed fiber $\mathscr L(T_{x_0}M,T_{y_0}N)$ and the function $f_{x_0}^{y_0}$ is continuous. Using also equation~\eqref{eqCarNec4}, we finally obtain
\begin{align*}
f(\alpha_{x_0}^{y_0})\leqslant&\,\frac{1}{\mu_g(Q_{x_0}^r)}\int_{Q_{x_0}^r}f\big(\alpha_{x_0}^{y_0}+d\exp_{y_0}^{-1}[\varphi(x)]\circ d\varphi[x]\circ d\exp_{x_0}[\exp_{x_0}^{-1}(x)]\big)\,d\mu_g(x)+o(1),
\end{align*}
for every $y_0\in N$, every $\alpha_{x_0}^{y_0}\in\mathscr L(T_{x_0}M,T_{y_0}N)$ and every admissible $\varphi$. The compact subsets in the argument above can be chosen increasing with respect to $\Vert d\varphi\Vert_\infty$ and the constants and moduli of continuity can be enlarged accordingly; hence the function $o(1)$ can be chosen with a monotonically nondecreasing dependence on $\Vert d\varphi\Vert_\infty$. Therefore, $f$ is R--quasiconvex at every $x_0\in M_0$, hence $\mu_g$--almost~everywhere in $M$.
\end{proof}

The necessity in the case $1\leqslant p<+\infty$ then follows easily.

\begin{proof}[Proof of the necessity -- The case $1\leqslant p<+\infty$]\ \\
Assume that, for every bounded open subset $\Omega$ of $M$, the functional $F(\cdot,\Omega)$ is sequentially lower semicontinuous with respect to the weak convergence in $W^{1,p}(\Omega,N)$. If
$$
u_j\weakkk u\qquad\text{in }W^{1,\infty}(\Omega,N),\qquad\text{ then }\qquad 
u_j\rightharpoonup u\qquad\text{in }W^{1,p}(\Omega,N).
$$
Here, as throughout the paper, $\rightharpoonup$ has the meaning fixed in Section~\ref{RDefSec}. Indeed, $\mu_g(\Omega)<+\infty$, as $\Omega$ is bounded and $M$ is complete; hence, if $1<p<+\infty$, the uniform convergence of $u_j$ to $u$ implies $d_{L^p}(u_j,u)\to0$, while the uniform boundedness of the differentials gives
$$
\sup_{j\in\NN}\Vert du_j\Vert_{L^p(\Omega)}<+\infty.
$$
If $p=1$, the same uniform convergence implies $d_{L^1}(u_j,u)\to0$ and the uniform
$L^\infty$ bound on $|du_j|$ implies the equiintegrability of the family
$|du_j|$. Hence, $F(\cdot,\Omega)$ is sequentially lower semicontinuous
also with respect to the weak$^*$ convergence in $W^{1,\infty}(\Omega,N)$. Moreover, if $K\subseteq N$ is compact and $R>0$, setting
$$
D_K=\max_{y\in K}d^N(y,y_*),
$$
the growth assumption gives
$$
0\leqslant f(\alpha_x^y)\leqslant a(x)+C\big(D_K^p+R^p\big)
$$
for $\mu_g$--almost~every $x\in M$, every $y\in K$ and every $\alpha_x^y$ with $|\alpha_x^y|\leqslant R$. Therefore, all the hypotheses needed in the necessity part of Theorem~\ref{AcerbiFuscoRinf} are satisfied and $f$ is R--quasiconvex.
\end{proof}

We now turn our attention to the sufficiency in Theorems~\ref{AcerbiFuscoRinf} and~\ref{AcerbiFuscoRp}.

\begin{proof}[Proof of the sufficiency -- The case $p=+\infty$]\ \\
Let $\Omega$ be a bounded open subset of $M$, and let us consider a sequence $u_j$ in $W^{1,\infty}(\Omega, N)$ and a map $u\in  W^{1,\infty}(\Omega, N)$ such that $u_j\weakkk u$ in $ W^{1,\infty}(\Omega, N)$. We aim to prove that
$$
F(u, \Omega)\leqslant \liminf_{j\to \infty} F(u_j, \Omega),
$$
so we can assume that the right--hand side is finite, since otherwise there is nothing to prove. Since $u_j\weakkk u$ in $ W^{1,\infty}(\Omega, N)$, by the definition of weak$^*$ convergence given in Section~\ref{RDefSec}, the sequence $u_j$ converges uniformly to $u$ and the norms $|du_j|$ are uniformly bounded in $L^{\infty}(\Omega)$. Thus, we may choose $R>0$ such that 
$$
|du[x]|\leqslant R\qquad\text{ and }\qquad |du_j[x]|\leqslant R,
$$
for $\mu_g$--almost~every $x\in\Omega$ and every sufficiently large $j\in\NN$. Since the image of $u$ is bounded and $N$ is complete, the Hopf--Rinow theorem and the uniform convergence of $u_j$ to $u$ then imply that there exists a compact subset $K\subseteq N$ containing the images of $u$ and of all $u_j$ for sufficiently large $j\in\NN$.\\
By the growth assumption~\eqref{growthAFinf}, we have $a\in L^1(\Omega)$ and there exists a constant $C_{K,R}\geqslant0$ such that
$$
0\leqslant f(du[x])\leqslant a(x)+C_{K,R}\qquad\text{ and }\qquad 0\leqslant f(du_j[x])\leqslant a(x)+C_{K,R}.
$$
Let $M_0\subseteq M$ be the set of full $\mu_g$--measure given by the Scorza~Dragoni argument in the proof of the necessity for $p=+\infty$: in every local ``trivialization'' $\mathscr L(TU,TV)$, for every $x_0\in M_0\cap U$ and every compact subset $\mathscr K\subseteq V\times\R^{m\times n}$ there exists a compact set $E\subseteq U$ containing $x_0$, of which $x_0$ is a density point, such that $f_{U,V}$ is continuous on $E\times\mathscr K$; moreover, every $x_0\in M_0$ is a Lebesgue point of $a$, the growth condition holds at $x_0$ and $f_{U,V}(x_0,\cdot,\cdot)$ is continuous. Since $f$ is R--quasiconvex, by possibly reducing $M_0$ by a set of measure zero, we can also assume that the R--quasiconvexity condition holds at every $x_0\in M_0$.\\
We fix a point $x_0\in M_0\cap\Omega$ which is also a Lebesgue point of the function $x\mapsto f(du[x])$ and we set $y_0=u(x_0)$. Moreover, we choose neighborhoods $U\subseteq\Omega$ of $x_0$ and $V\subseteq N$ of $y_0$, endowed with local orthonormal frames as in Section~\ref{quasisec}. Since $u\in W^{1,\infty}(\Omega,N)$, we can choose $r_0>0$ and a compact set $K'\subseteq V$ such that
$$
\overline{Q}_{x_0}^{r_0}\subseteq U\qquad\text{ and }\qquad u\big(\overline{Q}_{x_0}^{r_0}\big)\subseteq\mathop{\rm int}K'.
$$
Then, the uniform convergence of $u_j$ to $u$ implies that
$u_j\big(\overline{Q}_{x_0}^{r_0}\big)\subseteq \mathop{\rm int}K'$ for all sufficiently large $j\in\NN$. We set
$$
\mathscr K=K'\times\big\{A\in\R^{m\times n}\,\,\big|\,\,|A|\leqslant R\big\},
$$
and we let $E\subseteq U$ be a compact set associated with $\mathscr K$ by the Scorza~Dragoni construction above, such that $x_0\in E$ is a density point of $E$ and $f_{U,V}$ is continuous on $E\times\mathscr K$.\\
We now ``freeze'' the dependence of the integrand on the base point and consider on the local bundle $\mathscr L(TU,TV)$ the function
$$
g_{x_0}(\alpha_x^y)=f_{U,V}\big(x_0,y,I\alpha_x^y\big),
$$
which is continuous and R--quasiconvex: indeed, in the chosen local ``trivialization'', for every $x\in U$ and $y\in V$, its restriction to the fiber $\mathscr L(T_xM,T_yN)$ is represented by the quasiconvex function
$$
A\mapsto f_{U,V}(x_0,y,A),
$$
hence, the claim follows by the characterization of R--quasiconvexity in Proposition~\ref{RiemannianQC-equiv}. Notice that only the restriction of the integrand to this local bundle is involved in the argument below.\\
Thus, the local argument in the proof of Theorem~\ref{QCR->sci} applies without any modification to $g_{x_0}$ on $Q_{x_0}^r$, for every $0<r\leqslant r_0$ and we have
\begin{equation}\label{eqcar999}
\int_{Q_{x_0}^r}g_{x_0}(du[x])\,d\mu_g(x)\leqslant\liminf_{j\to\infty}\int_{Q_{x_0}^r}g_{x_0}(du_j[x])\,d\mu_g(x).
\end{equation}
Let now $\omega$ be a modulus of continuity for $f_{U,V}$ on the compact set $E\times\mathscr K$. Since the matrices representing $du$ and $du_j$ have norms bounded by $R$, for every $x\in E\cap Q_{x_0}^r$ and every sufficiently large $j\in\NN$, there holds
$$
\big|f(du[x])-g_{x_0}(du[x])\big|+\big|f(du_j[x])-g_{x_0}(du_j[x])\big|\leqslant2\omega(Cr),
$$
where $C$ is independent of $r$ and $j$. On the complementary set $Q_{x_0}^r\setminus E$, we use instead the growth assumption, both at $x$ and at the frozen point $x_0$, obtaining
$$
0\leqslant f(du[x])\leqslant a(x)+C_{K'\!,R}\qquad\text{ and }\qquad0\leqslant g_{x_0}(du_j[x])\leqslant a(x_0)+C_{K'\!,R}.
$$
Combining these estimates with the lower semicontinuity inequality~\eqref{eqcar999} for $g_{x_0}$ and using the nonnegativity of $f$ and $g_{x_0}$, we get
\begin{align*}
\int_{Q_{x_0}^r}f(du[x])\,d\mu_g(x)\leqslant&\,\liminf_{j\to\infty}\int_{Q_{x_0}^r}f(du_j[x])\,d\mu_g(x)+2\omega(Cr)\mu_g(Q_{x_0}^r)\\
&\,+\int_{Q_{x_0}^r\setminus E}\big[a(x)+C_{K'\!,R}\big]\,d\mu_g(x)+\big[a(x_0)+C_{K'\!,R}\big]\mu_g(Q_{x_0}^r\setminus E).
\end{align*}
Dividing both sides by $\mu_g(Q_{x_0}^r)$, the last two terms on the right tend to zero as $r\to0$, because $x_0$ is a density point of $E$ and a Lebesgue point of $a$. Since also $\omega(Cr)\to0$ and $x_0$ is a Lebesgue point of $f(du)$, it follows that
$$
f(du[x_0])\leqslant\liminf_{j\to\infty}\fint_{Q_{x_0}^r}f(du_j[x])\,d\mu_g(x)+o(1),
$$
as $r\to0$. Hence, we have obtained precisely the analogue of inequality~\eqref{eq558} for our Carath\'eodory integrand.\\
We can now conclude by repeating the measure--theoretic argument in the proof of Theorem~\ref{QCR->sci}, starting from this inequality; the only difference is that the uniform upper bound used there is replaced by the integrable bound
$$
0\leqslant f(du_j)\leqslant a+C_{K,R},
$$
indeed, after choosing a subsequence, not relabeled, such that
$$
\lim_{j\to\infty}F(u_j,\Omega)=\liminf_{j\to\infty}F(u_j,\Omega)
$$
and defining the finite positive Radon measures $d\nu_j=f(du_j)\,d\mu_g$, the preceding estimate gives a uniform bound on them. Up to a further subsequence, they converge weakly$^*$ to a finite positive Radon measure $\nu$. For every open set $G\subseteq\Omega$, we have
$$
\nu(G)\leqslant\liminf_{j\to\infty}\nu_j(G)
\leqslant\int_G\big(a(x)+C_{K,R}\big)\,d\mu_g(x).
$$
Hence, by the outer regularity of the finite Radon measure $(a+C_{K,R})\mu_g$, we have $\nu\ll\mu_g$. Writing
$d\nu=\theta\,d\mu_g$, the differentiation argument used in the proof of Theorem~\ref{QCR->sci} gives
$$
f(du[x])\leqslant\theta(x),
$$
for $\mu_g$--almost~every $x\in\Omega$. Therefore,
$$
F(u,\Omega)\leqslant\nu(\Omega)\leqslant\liminf_{j\to\infty}\nu_j(\Omega)=\liminf_{j\to\infty}F(u_j,\Omega),
$$
which proves the claim.
\end{proof}

Before proving the sufficiency for $p<+\infty$, we record some preliminary results. We fix an open set $\Omega'\subseteq M$ contained in a compact subset of $M$. On this compact region, the injectivity radius is bounded from below by a positive constant and, for sufficiently small radii, all the normal coordinate maps and their inverses have Lipschitz constants at most $2$; moreover, in these coordinates, the Riemannian volume measure has a density with respect to the Lebesgue measure bounded between $2^{-n}$ and $2^n$. Writing $\operatorname{inj}_M(x)$ for the injectivity radius of $M$ at $x$, we choose $r_0>0$ small enough for these uniform estimates to hold for every $0<r<r_0$ and so that
$$
\frac{\sqrt n}{2}\,r_0<\inf_{x\in\overline{\Omega}'}\operatorname{inj}_M(x).
$$

The following lemmas are respectively the Riemannian counterparts needed here of Lemmas~I.7, I.9 and I.11 in~\cite{AcerbiFusco}.

\begin{lem}[{\em Biting lemma}]\label{lemma1.7adattato}
Let $\varphi_k$ be a bounded sequence of nonnegative functions in $L^1(\Omega')$. Then, for every $\varepsilon>0$ there exists a triple $(A_\varepsilon, \delta, S)$, where $A_{\varepsilon}$ is measurable and $\mu_g(A_\varepsilon)<\varepsilon$, $\delta>0$ and $S$ is an infinite subset of $\NN$, such that for every $k\in S$, there holds
$$
\int_B \varphi_k(x)\,d\mu_g(x)<\varepsilon,
$$
whenever $B\subseteq \Omega'\setminus A_\varepsilon$ is measurable and $\mu_g(B)<\delta$.
\end{lem}

\begin{lem}\label{lemma1.9adattato}
Let $u\in W^{1,p}(\Omega',N)$ and let $y_*\in N$ be the fixed reference point in
Theorem~\ref{AcerbiFuscoRp}. We define
\begin{equation*}
h_u=d^N(u,y_*)+|du|
\end{equation*}
and 
\begin{equation*}
\mathcal{M}_g h_u(x)=\sup_{\substack{0<r<r_0\\Q_x^r\subseteq\Omega'}}
\fint_{Q_x^{r}}h_u\,d\mu_g,
\end{equation*}
where $r_0>0$ is the one fixed above.\\
Then, the following local estimates hold: if $p>1$, we have
$$
\|\mathcal{M}_gh_u \|_{L^p(\Omega')}\leqslant C(p,\Omega')\| h_u\|_{L^p(\Omega')}
$$
and, for every $1\leqslant p<+\infty$ and all $\lambda>0$,
$$
\mu_g\big(\big\{ x\in \Omega'\,\,\big|\,\,\mathcal{M}_gh_u(x) \geqslant \lambda\big \}\big)\leqslant\frac{C(\Omega')}{\lambda} \|h_u \|_{L^1(\Omega')}. 
$$
Moreover,
$$
h_u(x)\leqslant\mathcal{M}_gh_u(x)
$$
for $\mu_g$--almost~every $x\in\Omega'$.
\end{lem}

\begin{lem}\label{lemma1.11adattato}
Let $u\in W^{1,p}(\Omega',N)$, $\lambda>0$ and set
$$
H^\lambda=\big\{x\in\Omega'\,\,|\,\,\mathcal{M}_gh_u(x)<\lambda\big\}.
$$
Then, there exist constants $c_n>1$ and $C_1>0$, depending only on $n$, such that for every cube $Q=Q_{x_0}^r$ with $0<r<r_0/c_n$ and $\overline{Q}_{x_0}^{c_n r}\subseteq\Omega'$, we have
$$
d^N(u(x),u(z))\leqslant C_1\lambda d^M(x,z),
$$
for every $x,z\in H^\lambda\cap Q$ outside a fixed $\mu_g$--null set depending only on $u$.
\end{lem}

The proof of Lemma~\ref{lemma1.7adattato} is exactly the same as that of Lemma~I.7 of~\cite{AcerbiFusco}, with the Riemannian volume measure in place of the Lebesgue measure. 

About Lemma~\ref{lemma1.9adattato}, since by the choice of $r_0$ in the normal coordinates on the cubes $Q_x^r$, with $x\in\Omega'$ and $0<r<r_0$, we have uniform bounds on the Lipschitz constants of the coordinate maps and of their inverses and on the volume densities, as we said above, we can apply the Euclidean maximal--function estimates recalled in~\cite{AcerbiFusco}, obtaining the stated inequalities with constants depending only on $\Omega'$ and its geometry.

Regarding Lemma~\ref{lemma1.11adattato}, letting $D\subseteq N$ be a countable dense subset, for every $q\in D$ we consider the real--valued function
$$
v_q(x)=d^N(u(x),q)-d^N(y_*,q).
$$
Then, by the definition of $W^{1,p}(\Omega',N)$ in Section~\ref{RDefSec} and since $d^N(\cdot,q)$ is $1$--Lipschitz and $u$ is approximately differentiable almost everywhere, $v_q\in W^{1,p}(\Omega')$ and there holds
$$
|v_q|\leqslant d^N(u,y_*)\qquad\text{and}\qquad |dv_q|\leqslant |du|
$$
$\mu_g$--almost~everywhere in $\Omega'$. We now work in normal coordinates centered at $x_0$, in which $Q=Q_{x_0}^r$ and $Q_{x_0}^{c_nr}$ are the concentric Euclidean cubes of side $r$ and $c_nr$, and we regard $v_q$ as a function on the latter. The local $W^{1,p}$ version of Lemma~I.11 of~\cite{AcerbiFusco} (see~\cite{Liu}) on $Q_{x_0}^{c_nr}$ states that 
$$
|v_q(x)-v_q(z)|\leqslant c(n)\Lambda|x-z|
$$
whenever $x,z\in Q$ are Lebesgue points of $v_q$ at which the maximal function of $|v_q|+|Dv_q|$, taken over the Euclidean balls centered at the point with radius at most $2\sqrt n\,r$, is smaller than $\Lambda$; for $c_n$ large enough, depending only on $n$, these balls are contained in $Q_{x_0}^{c_nr}$. We claim that this maximal function is bounded by $C(n)\mathcal{M}_gh_u$ at every point $x\in Q$. Indeed, $|v_q|+|Dv_q|\leqslant2h_u$, by the two inequalities above on $v_q$ and $dv_q$ and since the normal coordinate maps are $2$--Lipschitz, as fixed before Lemma~\ref{lemma1.7adattato}; moreover, by the same Lipschitz bounds, every Euclidean ball centered at $x$ with radius $\rho'\leqslant2\sqrt n\,r$ is contained in the cube $Q_x^{8\rho'}$, which is contained in $Q_{x_0}^{c_nr}\subseteq\Omega'$ for $c_n$ large enough and satisfies $8\rho'<r_0$, since $r<r_0/c_n$; finally, the Lebesgue measure of such a ball and $\mu_g(Q_x^{8\rho'})$ are comparable, with constants depending only on $n$. Hence, the Euclidean averages of $|v_q|+|Dv_q|$ over these balls are bounded by $C(n)$ times the averages of $h_u$ over the cubes $Q_x^{8\rho'}$, which are among the cubes in the definition of $\mathcal{M}_gh_u(x)$. Applying the local version of Lemma~I.11 of~\cite{AcerbiFusco} stated above with $\Lambda=C(n)\lambda$ and using once more the Lipschitz bounds on the coordinates, we obtain
\begin{equation}\label{eq:lipvq}
|v_q(x)-v_q(z)|\leqslant C_1\lambda d^M(x,z)
\end{equation}
for every $q\in D$ and every $x,z\in H^\lambda\cap Q$ outside a fixed $\mu_g$--null set $Z_q$; moreover, $C_1=2c(n)C(n)$ depends only on $n$, in particular it is independent of $q$ and $Q$.\\
As $D$ is countable, the union $Z=\bigcup_{q\in D}Z_q$ is a $\mu_g$--null set which depends only on $u$, while inequality~\eqref{eq:lipvq} holds for every $q\in D$ whenever $x,z\in(H^\lambda\cap Q)\setminus Z$. Then, since
$$
d^N(y_1,y_2)=\sup_{q\in D}\big|d^N(y_1,q)-d^N(y_2,q)\big|,
$$
taking the supremum over $q\in D$, we conclude
$$
d^N(u(x),u(z))\leqslant C_1\lambda d^M(x,z)
$$
for every $x,z\in H^\lambda\cap Q$ outside $Z$ and this gives Lemma~\ref{lemma1.11adattato}.

\smallskip

Finally, before starting the proof of the sufficiency in Theorem~\ref{AcerbiFuscoRp}, we need one last preliminary result.

\begin{lem}\label{estensionelipschitz}
Let $K\subseteq N$ be a compact set. There exist $0<\rho_K\leqslant1$ and a constant $C_2\geqslant1$, depending only on $m$ and independent of $K$, with the following property. Let $Q=Q_{x_0}^r\subseteq\Omega'$ be a cube with $0<r<r_0$, $A\subseteq Q$, $L\geqslant0$ and $v:A\to K$ satisfy
$$
d^N(v(x), v(z))\leqslant Ld^M(x,z)
$$
for all $x,z\in A$. Then, if $L\diam Q\leqslant\rho_K$, there exists an extension $\widetilde{v}\in W^{1,\infty}(Q, N)$ such that
$$
\widetilde{v}=v \quad\text{on }A\qquad\text{ and }\qquad \| d\widetilde{v}\|_{L^{\infty}(Q)}\leqslant C_2 L.
$$
Moreover, if $A\neq\varempty$, then $\widetilde v(Q)$ is contained in the closed ball of radius $\rho_K$ centered at a point of $v(A)$.
\end{lem}
\begin{proof}
If $A=\varempty$, we choose any constant map on $Q$ and we are done. Otherwise, fix $a_0\in A$ and set $y_0=v(a_0)\in K$. By the compactness of $K$, after decreasing $\rho_K\leqslant1$ if necessary, we may assume that, for every $y\in K$, the Lipschitz constants of $\exp_y$ on the ball of radius $\rho_K$ centered at the origin of $T_yN$, where $\exp_y$ is injective and of its inverse are at most $2$ and notice that only the radius $\rho_K$ depends on $K$. If $L\diam Q\leqslant\rho_K$, then $v(A)$ is contained in the closed ball $\overline{B}_{\rho_K}(y_0)$, hence $\exp_{y_0}^{-1}\circ v$ is Lipschitz on $A$, with respect to $d^M$, with Lipschitz constant at most $2L$. Extending each of its components, in an orthonormal basis of $T_{y_0}N$, by the McShane extension theorem (see Lemma~I.12 in~\cite{AcerbiFusco}, which holds in every metric space), we obtain a map $w:Q\to T_{y_0}N$ with Lipschitz constant at most $2\sqrt m\,L$. Let $\pi:T_{y_0}N\to T_{y_0}N$ be the radial projection onto the closed ball of radius $\rho_K$ centered at the origin. Since $\pi$ fixes $\exp_{y_0}^{-1}(v(A))$ and does not increase the distances, the map
$$
\widetilde v=\exp_{y_0}\circ\pi\circ w:Q\to N
$$
extends $v$ from $A$ to $Q$ and, since $\exp_{y_0}$ is $2$--Lipschitz on the ball of radius $\rho_K$, we have $|d\widetilde v|\leqslant2|dw|\leqslant4\sqrt m\,L$
$\mu_g$--almost~everywhere on $Q$, using the $2L$--Lipschitz bound on each of the $m$ components of $w$, hence there holds
$$
\|d\widetilde v\|_{L^\infty(Q)}\leqslant C_2L,
$$
with $C_2=4\sqrt m$. Finally, $\widetilde v(Q)\subseteq\exp_{y_0}\big(\overline{B}_{\rho_K}(0)\big)\subseteq\overline{B}_{\rho_K}(y_0)$, with $y_0\in v(A)$.
\end{proof}

\begin{proof}[Proof of the sufficiency -- The case $1\leqslant p<+\infty$]\ \\
Let $\Omega'$ be a bounded open subset of $(M,g)$ and $u_k\in W^{1,p}(\Omega',N)$ be a sequence weakly converging to $u\in W^{1,p}(\Omega',N)$. Since $(M,g)$ is complete, by the Hopf--Rinow theorem, $\overline{\Omega}'$ is compact, hence the preceding lemmas apply. Possibly passing to a subsequence, we can assume that 
\begin{equation*}
\lim_{k\to \infty}F(u_k, \Omega')=\liminf_{k\to \infty} F(u_k, \Omega')<+\infty
\end{equation*}
and that $u_k$ converges $\mu_g$--almost everywhere to $u$. By the absolute continuity of the integral of the function $a$, fixed $\varepsilon>0$, we can choose $\delta>0$ such that
\begin{equation}\label{eq:small-a}
\int_B a\,d\mu_g<\varepsilon,
\end{equation}
whenever $\mu_g(B)<\delta$.
With the notation introduced in Lemma~\ref{lemma1.9adattato}, we set 
$$
h_k=h_{u_k}=d^N(u_k,y_*)+|du_k|,
$$
then by the same lemma, for every $\lambda>0$, we have
\begin{equation*}
\mu_g\big(\big\{x\in\Omega'\,\,\big|\,\,\mathcal{M}_g h_k(x)\geqslant\lambda\big\}\big)\leqslant \frac{C(\Omega')}{\lambda}\|h_k\|_{L^1(\Omega')}.
\end{equation*}
We want a uniform control on $ \| h_k\|_{L^1(\Omega')}$, for every $k\in\NN$: by means of the H\"older, Minkowski and triangle inequalities, we get
\begin{align*}
\| h_k\|_{L^1(\Omega')}=&\,\int_{{\Omega}'} \big(d^N(y_*,u_k(x))+|du_k[x]|\big)\,d\mu_g(x)\\
\leqslant&\,\mu_g(\Omega')^{1-\frac{1}{p}}[ \|d^N(y_*, u_k)\|_{L^p(\Omega')}+ \| du_k\|_{L^p(\Omega')}]\\
\leqslant&\,\mu_g(\Omega')^{1-\frac{1}{p}}[ \|d^N(y_*, u)\|_{L^p(\Omega')}+ \|d^N(u, u_k)\|_{L^p(\Omega')} +\| du_k\|_{L^p(\Omega')}],
\end{align*}
hence, there holds
$$
\| h_k\|_{L^1(\Omega')}\leqslant\sup_{k\in\NN}\big\{\mu_g(\Omega')^{1-\frac{1}{p}}[ \|d^N(y_*, u)\|_{L^p(\Omega')}+ \|d^N(u, u_k)\|_{L^p(\Omega')} +\| du_k\|_{L^p(\Omega')}] \big\}=\overline{C}<+\infty
$$
for every $k\in\NN$.\\
Now, for every $\lambda>0$, we set
$$
H_k^\lambda=\big\{x\in\Omega'\,\,|\,\,\mathcal M_g h_k(x)<\lambda\big\}.
$$
Our aim is to choose a measurable set $A_\varepsilon\subseteq\Omega'$ with $\mu_g(A_\varepsilon)<\varepsilon$ and a number $\lambda\geqslant1$, independent of $k$, such that
\begin{equation}\label{finale1}
\mu_g(\Omega'\setminus H_k^\lambda)<\varepsilon
\qquad\text{and}\qquad
\lambda^p\mu_g\big(\Omega'\setminus(H_k^\lambda\cup A_\varepsilon)\big)<\varepsilon,
\end{equation}
for every $k\in\NN$.\\
If $p>1$, the sequence $(\mathcal{M}_g h_k)^p$ is bounded in $L^1(\Omega')$, by Lemma~\ref{lemma1.9adattato}; then, applying Lemma~\ref{lemma1.7adattato} to this sequence, possibly passing to a further subsequence, we can suppose that there exist a measurable set $A_\varepsilon\subseteq\Omega'$ with $\mu_g(A_\varepsilon)<\varepsilon$ and a number $\delta'>0$ such that
\begin{equation}\label{eq:small-Psi}
\int_B(\mathcal{M}_g h_k)^p\,d\mu_g<\varepsilon,
\end{equation}
for every $k\in\NN$ and every measurable set $B\subseteq\Omega'\setminus A_\varepsilon$ with $\mu_g(B)<\delta'$. Replacing $\delta$ by $\min\{\delta,\delta'\}$, inequalities~\eqref{eq:small-a} and~\eqref{eq:small-Psi} hold whenever $B\subseteq\Omega'\setminus A_\varepsilon$ is measurable and $\mu_g(B)<\delta$.\\
Moreover, the $L^p$ bound on $\mathcal M_g h_k$ also gives, by the Chebyshev inequality,
$$
\sup_{k\in\NN}\mu_g(\Omega'\setminus H_k^\lambda)\leqslant\frac{1}{\lambda^p}
\sup_{k\in\NN}\int_{\Omega'}(\mathcal M_g h_k)^p\,d\mu_g\to0,
$$
as $\lambda\to+\infty$. Thus, for every sufficiently large $\lambda>0$, we have
$$
\mu_g(\Omega'\setminus H_k^\lambda)<\min\{\varepsilon,\delta\}
$$
for every $k\in\NN$ and we can therefore apply inequality~\eqref{eq:small-Psi} to the measurable set $\Omega'\setminus (H_k^\lambda\cup A_\varepsilon)$, whose measure is smaller than $\delta$. Since $\mathcal M_g h_k\geqslant\lambda$ on $\Omega'\setminus H_k^\lambda$, we obtain
$$
\lambda^p\mu_g\big(\Omega'\setminus (H_k^\lambda\cup A_\varepsilon)\big)
\leqslant\int_{\Omega'\setminus (H_k^\lambda\cup A_\varepsilon)}
(\mathcal M_g h_k)^p\,d\mu_g
<\varepsilon.
$$
Together with the bound $\mu_g(\Omega'\setminus H_k^\lambda)<\varepsilon$, we have both inequalities~\eqref{finale1} for $p\in(1,+\infty)$.\\
The case $p=1$ requires a different treatment, since the maximal operator is not bounded from $L^1(\Omega')$ to itself, hence, we cannot deduce the $L^1$ boundedness of $\mathcal M_g h_k$ needed in order to apply Lemma~\ref{lemma1.7adattato}.\\
We set $A_\varepsilon=\varempty$ and we notice that the sequence $h_k$ is equiintegrable, indeed,
$$
h_k\leqslant d^N(u_k,u)+d^N(u,y_*)+|du_k|,
$$
where the sequence $d^N(u_k,u)$ converges to zero in $L^1(\Omega')$, the function $d^N(u,y_*)$ belongs to $L^1(\Omega')$ and the sequence $|du_k|$ is equiintegrable by our notion of weak convergence.\\
We then apply the argument leading to formula~(3.10) in~\cite{FocardiMarcellini}. For every $k\in\NN$ and $\lambda>0$, we set
$$
h_{k,\lambda}=h_k\rchi_{\{h_k>\lambda/2\}}
$$
and we define its maximal function by
$$
\mathcal M_g h_{k,\lambda}(x)=\sup_{\substack{0<r<r_0\\Q_x^r\subseteq\Omega'}}
\fint_{Q_x^r}h_{k,\lambda}\,d\mu_g
$$
for every $x\in\Omega'$. Since $h_{k,\lambda}$ is a nonnegative function in $L^1(\Omega')$, we have, as in Lemma~\ref{lemma1.9adattato},
$$
\mu_g\big(\{x\in\Omega'\,\,|\,\,\mathcal M_g h_{k,\lambda}(x)\geqslant\lambda/2\}\big)
\leqslant\frac{2C(\Omega')}{\lambda}\int_{\Omega'}h_{k,\lambda}\,d\mu_g
=\frac{2C(\Omega')}{\lambda}\int_{\{h_k>\lambda/2\}}h_k\,d\mu_g
$$
and the bound $\|h_k\|_{L^1(\Omega')}\leqslant \overline{C}$ implies 
$$
\mu_g(\{h_k>\lambda/2\})\leqslant2\overline{C}/\lambda.
$$
Hence, by the equiintegrability of the sequence $h_k$, the last integrals on the right in the above formula tend to zero uniformly in $k$, as $\lambda\to+\infty$.\\
Since $0\leqslant h_k-h_{k,\lambda}\leqslant\lambda/2$, we have
$$
\mathcal M_g h_k\leqslant\frac{\lambda}{2}+\mathcal M_g h_{k,\lambda},
$$
hence,
$$
\{x\in\Omega'\,\,|\,\,\mathcal M_g h_k(x)\geqslant\lambda\}
\subseteq
\{x\in\Omega'\,\,|\,\,\mathcal M_g h_{k,\lambda}(x)\geqslant\lambda/2\},
$$
we then conclude that
\begin{equation}\label{eq:maximal-p1}
\sup_{k\in\NN}\lambda\,\mu_g\big(\{x\in\Omega'\,\,|\,\,\mathcal M_g h_k(x)\geqslant\lambda\}\big)\to0,
\end{equation}
as $\lambda\to+\infty$.\\
By limit~\eqref{eq:maximal-p1}, we can choose $\lambda\geqslant1$, independently of $k$, such that
$$
\lambda\mu_g(\Omega'\setminus H_k^\lambda)<\min\{\varepsilon,\delta\},
$$
for every $k\in\NN$ and since $A_\varepsilon=\varempty$, this gives the second inequality~\eqref{finale1} for $p=1$. As $\lambda\geqslant1$, we have
$$
\mu_g(\Omega'\setminus H_k^\lambda)\leqslant\lambda\mu_g(\Omega'\setminus H_k^\lambda)
<\min\{\varepsilon,\delta\},
$$
hence, the first inequality~\eqref{finale1} holds as well.\\
Thus, in both cases, that is, for every $1\leqslant p<+\infty$, we have a measurable set $A_\varepsilon$ and a constant $\lambda\geqslant1$ so that inequalities~\eqref{finale1} hold and $\mu_g(\Omega'\setminus H_k^\lambda)<\delta$ also, for every $k\in\NN$. Consequently, inequality~\eqref{eq:small-a} gives
\begin{equation}\label{stima3}
\int_{\Omega'\setminus (H_k^\lambda\cup A_\varepsilon)}a\,d\mu_g<\varepsilon,
\end{equation}
for every $k\in\NN$.\\
Let $Z_k$ be the union of the $\mu_g$--null sets given by Lemmas~\ref{lemma1.9adattato} and~\ref{lemma1.11adattato} applied to $u_k$. The union $Z=\bigcup_{k\in\NN}Z_k$ is a $\mu_g$--null set and we replace each $H_k^\lambda$ by $H_k^\lambda\setminus Z$, keeping the same notation; then,
$$
d^N(y_*,u_k(x))\leqslant h_k(x)\leqslant\mathcal M_g h_k(x)<\lambda,
$$
for every $k\in\NN$ and every $x\in H_k^\lambda$, hence,
$$
u_k(H_k^\lambda)\subseteq\overline{B}_\lambda(y_*),
$$
where $\overline{B}_\lambda(y_*)$ is the closed ball in $N$ with center $y_*$ and radius $\lambda$ (which is compact by the Hopf--Rinow theorem).\\
Let $0<\rho_\lambda\leqslant1$ be the radius given by Lemma~\ref{estensionelipschitz} for the compact set $\overline{B}_\lambda(y_*)$; the constant $C_2$ in the extension estimate is independent of $\lambda$. To apply Lemma~\ref{lemma1.11adattato} and then Lemma~\ref{estensionelipschitz} with $L=C_1\lambda$, we consider the family
$$
\mathcal V=\big\{Q=Q_x^r\,\,\big|\,\,\text{$x\in\Omega'$, $0<r<r_0/c_n$, $\overline{Q}_{x}^{c_n r}\subseteq\Omega'$, $C_1\lambda\diam Q\leqslant\rho_\lambda$}\big\}
$$
which is clearly a Vitali covering of $\Omega'$, hence, by the Vitali covering theorem, we can choose a countable family of pairwise disjoint cubes $Q_i\in\mathcal V$ such that
$$
\mu_g\bigg(\Omega'\setminus\bigcup_{i\in\NN}Q_i\bigg)=0.
$$
For every $i\in\NN$, by Lemma~\ref{lemma1.11adattato}, for every $x,z\in H_k^\lambda\cap Q_i$, we have
$$
d^N(u_k(x),u_k(z))\leqslant C_1\lambda d^M(x,z),
$$
then, if $H_k^\lambda\cap Q_i=\varempty$, we simply set $w_{k,i}= y_*$ on $Q_i$, otherwise, applying Lemma~\ref{estensionelipschitz} on the cube $Q_i$ to
$$
A=H_k^\lambda\cap Q_i,\qquad L=C_1\lambda\qquad\text{ and }\qquad K=\overline{B}_\lambda(y_*),
$$
we obtain an extension $w_{k,i}\in W^{1,\infty}(Q_i,N)$ such that
$$
w_{k,i}=u_k \quad\text{on }H_k^\lambda\cap Q_i,
$$
for $\mu_g$--almost~every point and
\begin{equation}\label{stima1}
\|dw_{k,i}\|_{L^\infty(Q_i)}\leqslant C_1C_2\lambda.
\end{equation}
If $H_k^\lambda\cap Q_i\neq\varempty$, by Lemma~\ref{estensionelipschitz}, the set $w_{k,i}(Q_i)$ is contained in the closed ball of radius $\rho_\lambda \leqslant1$ centered at a point of $\overline{B}_\lambda(y_*)$, while if $H_k^\lambda\cap Q_i=\varempty$ we have $w_{k,i}= y_*$. Since $\lambda\geqslant1$, in either case, there holds
\begin{equation}\label{stima2}
d^N(y_*,w_{k,i})\leqslant\lambda+1\leqslant2\lambda
\end{equation}
on $Q_i$.\\
The sequences $w_{k,i}$, for every fixed $i\in\NN$, have uniformly bounded differentials on $Q_i$ and their images are contained in the compact set $\overline B_{2\lambda}(y_*)$, hence, for every $i\in\NN$, the maps $w_{k,i}$ extend continuously to $\overline{Q}_i$ with Lipschitz constants bounded independently of $k$ and values in the compact set $\overline{B}_{2\lambda}(y_*)$, so the Ascoli--Arzel\`a theorem gives a subsequence converging uniformly to a Lipschitz map $v_i\in W^{1,\infty}(Q_i,N)$. Then, by a diagonal argument, we obtain a subsequence, not relabeled, such that $w_{k,i}$ converges uniformly to maps $v_i\in W^{1,\infty}(Q_i,N)$ on $Q_i$, as $k\to\infty$, for every $i\in\NN$. Moreover, since the differentials $dw_{k,i}$ are uniformly bounded, we can also assume that $w_{k,i}\weakkk v_i$ in $W^{1,\infty}(Q_i,N)$.\\
Now, since $w_{k,i}=u_k$ pointwise $\mu_g$--almost~everywhere on $H_k^\lambda\cap Q_i$, we have that $\mu_g$--almost~everywhere on the same set, $dw_{k,i}=du_k$, hence, for every $i\in\NN$, there holds
\begin{equation}\label{liminf1}
F(u_k, Q_i)\geqslant F(w_{k,i},(Q_i\setminus A_{\varepsilon})\cap H^{\lambda}_k )=F(w_{k,i}, Q_i\setminus A_{\varepsilon})- F(w_{k,i},Q_i\setminus (H^{\lambda}_k\cup A_{\varepsilon})).
\end{equation}
For every $i,k\in\NN$ and $\ell\in\NN$, we set
$$
B_{k,i}=Q_i\setminus (H_k^\lambda\cup A_{\varepsilon})\qquad\text{ and }\qquad B_k^{\ell}=\bigcup_{i=1}^\ell B_{k,i}.
$$
Since the cubes $Q_i$ are pairwise disjoint, the sets $B_{k,i}$ are pairwise disjoint as well and $B_k^{\ell}\subseteq\Omega'\setminus (H^{\lambda}_k\cup A_{\varepsilon})$. By the growth condition on $f$ (inequality~\eqref{growthAFp}) and estimates~\eqref{stima1} and~\eqref{stima2}, we have
\begin{equation}\label{bound_f_wki}
f(dw_{k,i}[x])\leqslant a(x)+C\big(d^N(y_*,w_{k,i}(x))^p+|dw_{k,i}[x]|^p\big)\leqslant a(x)+C_3\lambda^p
\end{equation}
for $\mu_g$--almost~every $x\in Q_i$ and every $i,k\in\NN$, where $C_3=C(2^p+C_1^pC_2^p)$ is independent of $k$, $i$ and $\varepsilon$. Integrating this bound over $B_k^{\ell}$ and using inequalities~\eqref{stima3} and~\eqref{finale1}, we obtain
\begin{equation}\label{stima5}
\sum_{i=1}^\ell F(w_{k,i}, B_{k,i})\leqslant\int_{\Omega'\setminus (H^{\lambda}_k\cup A_{\varepsilon})}a\,d\mu_g+C_3\lambda^p\mu_g\big(\Omega'\setminus (H^{\lambda}_k\cup A_{\varepsilon})\big)\leqslant (1+C_3)\varepsilon,
\end{equation}
for every $\ell,k\in\NN$.\\
Moreover, for fixed $\varepsilon$, by estimate~\eqref{bound_f_wki} and the inner regularity of the finite Radon measure $(a+C_3\lambda^p)\mu_g$, we can choose a compact set $K_{i,\varepsilon}\subseteq Q_i\cap A_\varepsilon$ such that
$$
\int_{(Q_i\cap A_\varepsilon)\setminus K_{i,\varepsilon}}(a+C_3\lambda^p)\,d\mu_g<\varepsilon 2^{-i}.
$$
Setting $Q_i''=Q_i\setminus K_{i,\varepsilon}$ (notice that the sets $Q_i''$ are open), we have $Q_i\setminus A_\varepsilon\subseteq Q_i''\subseteq Q_i$, hence, for every $k\in\NN$,
\begin{equation}\label{stima6}
0\leqslant F(w_{k,i},Q_i'')-F(w_{k,i},Q_i\setminus A_\varepsilon)=F(w_{k,i},(Q_i\cap A_\varepsilon)\setminus K_{i,\varepsilon})<\varepsilon 2^{-i}.
\end{equation}
Applying Theorem~\ref{AcerbiFuscoRinf} (whose growth assumption follows from the present one as in the proof of the necessity) to the sequence $w_{k,i}\weakkk v_i$ in $W^{1,\infty}(Q_i'',N)$, we have
\begin{equation}\label{lsc_Qi}
\liminf_{k\to\infty}F(w_{k,i},Q_i'')\geqslant F(v_i,Q_i''),
\end{equation}
for every $i\in\NN$. Then, by inequalities~\eqref{liminf1},~\eqref{stima5} and~\eqref{stima6}, we obtain
\begin{align*}
F(u_k,\Omega')\geqslant&\,\sum_{i=1}^{\ell}F(u_k,Q_i)\geqslant \sum_{i=1}^{\ell}F(w_{k,i},Q_i'')-\varepsilon\sum_{i=1}^{\ell}2^{-i}-\sum_{i=1}^{\ell}F(w_{k,i},B_{k,i})\\
\geqslant&\,\sum_{i=1}^{\ell}F(w_{k,i},Q_i'')-(2+C_3)\varepsilon.
\end{align*}
Passing to the $\liminf$ as $k\to\infty$ and using estimate~\eqref{lsc_Qi}, we deduce
\begin{equation}\label{stima_liminf_uk}
\liminf_{k\to\infty}F(u_k,\Omega')\geqslant\sum_{i=1}^{\ell}F(v_i,Q_i'')-(2+C_3)\varepsilon.
\end{equation}
For every $i\in\NN$, we now define the two sets
$$
P_i=\{ x\in Q_i\,\,|\,\,v_i(x)\neq u(x) \}\qquad\text{ and }\qquad\widetilde{P}_i=P_i\cap \{x\in Q_i\,\,|\,\,u_k(x)\to u(x)\},
$$
which have the same measure, as $u_k$ converges $\mu_g$--almost~everywhere to $u$. Moreover, since $w_{k,i}$ converges uniformly to $v_i$ on $Q_i$ and $w_{k,i}=u_k$ $\mu_g$--almost~everywhere on $H_k^\lambda\cap Q_i$, we have, up to a fixed $\mu_g$--null set,
$$
\widetilde{P}_i\subseteq\liminf_{k\to\infty}(Q_i\setminus H_k^\lambda),
$$
hence,
$$
\mu_g(P_i)=\mu_g(\widetilde{P}_i)\leqslant\liminf_{k\to\infty}\mu_g(Q_i\setminus H_k^\lambda).
$$
Then, setting
$$
G=\bigcup_{i\in\NN}P_i,
$$
by Fatou's lemma, the disjointness of the cubes $Q_i$ and the first inequality in~\eqref{finale1}, there holds
$$
\mu_g(G)=\sum_{i\in\NN}\mu_g(P_i)\leqslant\liminf_{k\to\infty}\sum_{i\in\NN}\mu_g(Q_i\setminus H_k^\lambda)\leqslant\varepsilon,
$$
thus, $\mu_g(A_\varepsilon\cup G)\leqslant2\varepsilon$ and the absolute continuity of the integral of $f(du)$ gives a function $\sigma(\varepsilon)$, with $\sigma(\varepsilon)\to0$ as $\varepsilon\to0$, such that
\begin{equation}\label{stimaG}
F(u,A_\varepsilon\cup G)\leqslant\sigma(\varepsilon).
\end{equation}
Now, since $v_i=u$ pointwise on $Q_i\setminus P_i$, we have that $\mu_g$--almost~everywhere on the same set, $dv_i=du$. As $Q_i\setminus Q_i''\subseteq Q_i\cap A_\varepsilon$, we have
$$
Q_i\setminus (Q_i''\setminus P_i)=(Q_i\setminus Q_i'')\cup (Q_i''\cap P_i)\subseteq (Q_i\cap A_\varepsilon)\cup P_i,
$$
therefore,
$$
F(v_i, Q_i'')\geqslant F(v_i, Q_i''\setminus P_i)=F(u, Q_i''\setminus P_i)\geqslant F(u, Q_i)-F\big(u,(Q_i\cap A_\varepsilon)\cup P_i\big).
$$
Since the sets $(Q_i\cap A_\varepsilon)\cup P_i\subseteq Q_i$ are pairwise disjoint and contained in $A_\varepsilon\cup G$, summing this inequality over $i=1,\dots,\ell$ and using estimate~\eqref{stimaG}, we get
\begin{equation}\label{stima_somma_vi}
\sum_{i=1}^\ell F(v_i,Q_i'')\geqslant\sum_{i=1}^\ell F(u,Q_i)-F(u,A_\varepsilon\cup G)\geqslant\sum_{i=1}^\ell F(u,Q_i)-\sigma(\varepsilon).
\end{equation}
Combining inequalities~\eqref{stima_liminf_uk} and~\eqref{stima_somma_vi}, for every $\ell\in\NN$ we obtain
$$
\liminf_{k\to\infty}F(u_k,\Omega')\geqslant\sum_{i=1}^\ell F(u,Q_i)-\sigma(\varepsilon)-(2+C_3)\varepsilon.
$$
Letting $\ell\to\infty$ and using the fact that
$$
\mu_g\bigg(\Omega'\setminus\bigcup_{i\in\NN}Q_i\bigg)=0,
$$
together with the absolute continuity of $F(u,\cdot)$ with respect to $\mu_g$, we have
$$
\liminf_{k\to\infty}F(u_k,\Omega')\geqslant F(u,\Omega')-\sigma(\varepsilon)-(2+C_3)\varepsilon
$$
and finally sending $\varepsilon\to0$, we conclude
$$
\liminf_{k\to\infty}F(u_k,\Omega')\geqslant F(u,\Omega'),
$$
which proves the weak sequential lower semicontinuity of $F$ in $W^{1,p}(\Omega',N)$.
\end{proof}

\begin{oss}\label{FSpone}
For $1\leqslant p\leqslant+\infty$, the notion of convergence considered by Focardi and Spadaro~\cite{FocardiSpadaro} for maps $u_k,u\in W^{1,p}(\Omega,N)$ is given by
\begin{equation}\label{eq:FSconvergence}
\Vert d^N(u_k,u)\Vert_{L^p(\Omega)}\to0
\qquad\text{and}\qquad
\sup_{k\in\NN}\Vert du_k\Vert_{L^p(\Omega)}<+\infty.
\end{equation}
As discussed in Remark~\ref{convFS}, for $1<p\leqslant+\infty$ this convergence coincides with ours, defined in Section~\ref{RDefSec}. At $p=1$, instead, the convergence defined by conditions~\eqref{eq:FSconvergence} is weaker than ours, since it does not require equiintegrability of the differentials.\\
This assumption cannot be dropped from the statement in Theorem~\ref{AcerbiFuscoRp} for Carath\'eodory integrands. Indeed, if $M=N=\R$ and $\Omega=(0,1)$, letting $K\subseteq(0,1)$ be a compact set with empty interior and positive Lebesgue measure and setting $A=(0,1)\setminus K$, we consider the integrand
$$
f(x,s,\xi)=\big(1+\rchi_K(x)\big)|\xi|,
$$
which is a Carath\'eodory function (clearly not continuous), is convex in $\xi$ and satisfies
$$
|\xi|\leqslant f(x,s,\xi)\leqslant2|\xi|,
$$
with the associated functional
$$
u\mapsto F(u,(0,1))=\int_{(0,1)}\big(1+\rchi_K\big)|u'|\,d\mathscr{L}^1
$$
defined for every $u\in W^{1,1}(0,1)$.\\
Then, for every positive integer $k$ and $j=0,\ldots,k-1$, we define
$$
I_{k,j}=\big(j/k,(j+1)/k\big)\qquad\text{ and }\qquad\ell_{k,j}=\mathscr{L}^1(A\cap I_{k,j})>0
$$
(the inequality holds because $A$ is open and dense) and we set
$$
v_k=\sum_{j=0}^{k-1}\frac{\rchi_{A\cap I_{k,j}}}{k\ell_{k,j}}\qquad\text{ and }\qquad u_k(x)=\int_{(0,x)}v_k\,d\mathscr{L}^1,
$$
for every $x\in(0,1)$. Thus, we have $u_k\in W^{1,\infty}(0,1)$ and $u_k'=v_k\geqslant0$ almost everywhere. Identifying each $u_k$ with its continuous extension to $[0,1]$, we also have
$$
\int_{I_{k,j}}v_k\,d\mathscr{L}^1={1}/{k}\qquad\text{and}\qquad
u_k(j/k)=j/k,
$$
for every $j=0,\ldots,k-1$, with $u_k(1)=1$.\\
Hence, letting $u(x)=x$, it is easy to check that by construction (in particular, by the monotonicity of $u_k$),
$$
\Vert u_k-u\Vert_{L^\infty(0,1)}\leqslant{1}/{k}\qquad\text{and}\qquad\Vert u_k'\Vert_{L^1(0,1)}=1,
$$
that is, conditions~\eqref{eq:FSconvergence} hold for $p=1$; moreover, notice that every $u_k$ has the same boundary values as $u$, namely zero at $0$ and one at $1$. Since $v_k=0$ almost everywhere on $K$, we conclude
$$
F(u_k,(0,1))=\int_{(0,1)}\big(1+\rchi_K\big)v_k\,d\mathscr{L}^1=\int_{(0,1)}v_k\,d\mathscr{L}^1=1
$$
while
$$
F(u,(0,1))=\int_{(0,1)}\big(1+\rchi_K\big)\,d\mathscr{L}^1=1+\mathscr{L}^1(K)>1,
$$
thus, the functional $F$ is not sequentially lower semicontinuous with respect to the convergence defined by conditions~\eqref{eq:FSconvergence} for $p=1$.\\
We remark that the sequence $v_k$ is not equiintegrable. Indeed, $v_k\to\rchi_A$ almost everywhere: on $K$ all the $v_k$ vanish, while for almost every $x\in A$ the interval $I_{k,j}$ containing $x$ is contained in $A$ for $k$ large enough, hence $v_k(x)=1$. By Vitali's theorem, equiintegrability would imply convergence in $L^1(0,1)$; therefore one would get
$$
1=\lim_{k\to\infty}\int_{(0,1)}v_k\,d\mathscr{L}^1=\mathscr{L}^1(A)=1-\mathscr{L}^1(K)<1,
$$
which is a contradiction.
\end{oss}

\section{Conclusion and future directions}\label{conclusion}

We have introduced an intrinsic notion of R--quasiconvexity for integrands on $\mathscr L(TM,TN)$ and shown that it provides the natural Riemannian counterpart of Morrey's quasiconvexity for maps between smooth Riemannian manifolds. For continuous integrands, R--quasiconvexity characterizes sequential weak$^*$ lower semicontinuity in $W^{1,\infty}(\Omega,N)$. Under the corresponding Riemannian growth assumptions, this characterization extends to Carath\'eodory integrands and to the full range $1\leqslant p\leqslant+\infty$, yielding Riemannian analogues of the classical Acerbi--Fusco theorems. For $p=1$, the convergence considered in Theorem~\ref{AcerbiFuscoRp} requires
the equiintegrability of the differentials, a condition which cannot in general be dropped for Carath\'eodory integrands, as shown in Remark~\ref{FSpone}. The formulation also clarifies how the intrinsic convergence of manifold--valued Sobolev maps interacts with the local linear structure supplied by tangent spaces and exponential coordinates.

These results suggest several directions for further investigation. A first direction is to understand how much of the theory survives for nonsmooth or metric targets. The Sobolev and convergence notions used here are already largely metric in nature, in the spirit of the theories developed in~\cite{Chiron,KorevaarSchoen,Reshetnyak}, but R--quasiconvexity itself relies essentially on tangent spaces and exponential maps. An extension to metric targets would therefore require an appropriate replacement for this infinitesimal linear structure, possibly through metric differentials or other notions of tangent object.

A second direction concerns the relation between the intrinsic formulation developed here and extrinsic approaches based on isometric embeddings, tangential quasiconvexity and relaxation, as in~\cite{DFMT}. A more systematic comparison could lead to representation and relaxation results for manifold--constrained energies and may clarify which parts of the theory are genuinely intrinsic and which are invariant reformulations of Euclidean constrained problems.

It would also be natural to investigate higher--order functionals and weaker regularity classes, in particular Riemannian analogues of the classical higher--order theory~\cite{Meyers} and of relaxation in $BV$--type spaces for manifold--valued maps, in connection with the metric theories of~\cite{AmbrosioBV,Chiron} and the constrained frameworks recalled above.

The lower semicontinuity results proved here provide one of the main ingredients of the direct method. Combined with suitable coercivity and compactness assumptions, they suggest further questions concerning existence, relaxation and regularity of minimizers for anisotropic and non--quadratic variational problems on curved domains with curved targets; in this direction, the polyconvex existence theory for Cartesian maps developed in~\cite{BB} provides a natural point of comparison.

\bibliographystyle{amsplain}
\bibliography{References}

\end{document}